\documentclass[12pt]{article}
\usepackage{latexsym,amsmath}
\usepackage{amssymb}
\usepackage{bm}
\usepackage{subfigure}
\usepackage{graphicx}
\usepackage{cite}
\usepackage{algorithm,algorithmic}
\usepackage{url}
\usepackage{booktabs}
\date{}

\newtheorem{remark}{Remark}[section]

\newtheorem{theorem}{Theorem}[section]
\newtheorem{lemma}{Lemma}[section]

\numberwithin{equation}{section}
\begin{document}
\title{Fast iterative method with a second order implicit difference scheme for
time-space fractional convection-diffusion equations}
\author{Xian-Ming Gu$^{1,2}$\thanks{E-mail: guxianming@live.cn, x.m.gu@rug.nl.},~
Ting-Zhu Huang$^{1}$\thanks{Corresponding author. E-mail: tingzhuhuang@126.com.
Tel.: +86 28 61831608.},~Cui-Cui Ji$^{3}$\thanks{E-mail: cuicuiahuan@163.com.},\\
Bruno Carpentieri$^{4}$\thanks{E-mail: bcarpentieri@gmail.com.},~Anatoly A.
Alikhanov$^{5}$\thanks{E-mail: aaalikhanov@gmail.com.}
\\
{\small{\it 1. School of Mathematical Sciences,}}\\
{\small{\it University of Electronic Science and Technology of China,}}\\
{\small{\it Chengdu, Sichuan 611731, P.R. China}}\\
{\small{\it 2. Institute of Mathematics and Computing Science, University of
Groningen,}} \\
{\small{\it Nijenborgh 9, P.O. Box 407, 9700 AK Groningen, The Netherlands}}\\
{\small{\it 3. Department of Mathematics,}}\\
{\small{\it Southeast University, Nanjing, Jiangsu 210096, P.R. China}}\\
{\small{\it 4. School of Science and Technology,}}\\
{\small{\it Nottingham Trent University, Clifton Campus, Nottingham, NG11 8NS,
UK}}\\
{\small{\it 5. Institute of Applied Mathematics and Automation,}}\\
{\small{\it Russian Academy of Sciences, ul. Shortanova 89 a, Nalchik, 360000,
Russia}}
}

\maketitle

\begin{abstract}
In this paper we want to propose practical numerical methods to solve a class of
initial-boundary problem of time-space fractional convection-diffusion equations
(TSFCDEs). To start with, an implicit difference method based on two-sided weighted
shifted Gr\"{u}nwald formulae is proposed with a discussion of the stability and
convergence. We construct an implicit difference scheme (IDS) and show that it
converges with second order accuracy in both time and space. Then, we develop fast
solution methods for handling the resulting system of linear equation with the
Toeplitz matrix. The fast Krylov subspace solvers with suitable circulant preconditioners
are designed to deal with the resulting Toeplitz linear systems. Each time level
of these methods reduces the memory requirement of the proposed implicit difference scheme
from $\mathcal{O}(N^2)$ to $\mathcal{O}(N)$ and the computational complexity from
$O(N^3)$ to $O(N\log N)$ in each iterative step, where $N$ is the number of grid nodes. Extensive numerical example runs
show the utility of these methods over the traditional direct solvers of the implicit
difference methods, in terms of computational cost and memory requirements.


\emph{AMSC (2010)}: 65F15, 65H18, 15A51.
\end{abstract}

\section{Introduction}
\label{intro}
\quad\
In recent years there has been a growing interest in the field of fractional calculus.
Podlubny \cite{IPFDE}, Samko et al. \cite{SSAKO} and Kilbas et al. \cite{AAKHMS} provide
the history and a comprehensive treatment of this subject. Many phenomena in engineering,
physics, chemistry and other sciences can be described very successfully by using fractional
differential equations (FDEs). Diffusion with an additional velocity field and diffusion
under the influence of a constant external force field are, in the Brownian case, both modelled by the
convection-dispersion equation. In the case of anomalous diffusion this is no longer true, i.e.,
the fractional generalisation may be different for the advection case and the transport an external
force field in \cite{RMJKT}. In the study, we are very interested in the fast solver for solving the
initial-boundary value problem of the time-space fractional convection-diffusion
equation (TSFCDE) \cite{AISGMZ,YPST}:
\begin{equation}
\begin{cases}
\partial^{\alpha}_{0,t}u(x, t) = \gamma(t)\frac{\partial u(x,t)}{\partial x}
+ d_{+}(t){}_aD^{\beta}_{x}u(x,t) + d_{-}(t){}_xD^{\beta}_{b}u(x,t) + f(x,t),\\
u(x, 0) = \phi(x),\qquad\qquad\quad a \leq x\leq b,\\
u(a, t) = u(b, t) = 0,\qquad\quad 0 < t\leq T,\\
\end{cases}
\label{eq:lin}
\end{equation}
where $\alpha\in (0, 1],~\beta\in (1, 2],~a < x < b$, and $0 < t \leq T$. Here,
the parameters $\alpha$ and $\beta$ are the order of the TSFCDE, $f(x, t)$ is the
source term, and  diffusion coefficient functions $d_{\pm}(t)$ are non-negative
under the assumption that the flow is from left to right. Moreover, the variable
coefficients $\gamma(t)$ are real. The TSFCDE (\ref{eq:lin}) can be regarded as generalizations
of classical convection-diffusion equations with the first-order time derivative
replaced by the Caputo fractional derivative of order $\alpha\in (0, 1]$, and the
second-order space derivatives replaced by the two-sided Riemman-Liouville fractional
derivative of order $\beta\in (1, 2]$. Namely, the time fractional derivative in (\ref{eq:lin})
is the Caputo fractional derivative of order $\alpha$ \cite{IPFDE} denoted by
\begin{equation}
\partial^{\alpha}_{0,t}u(x, t) = \frac{1}{\Gamma(1 - \alpha)}\int^{t}_0 \frac{\partial
u(x,\xi)}{\partial \xi}\frac{d\xi}{(t - \xi)^{\alpha}},
\label{eq1.2}
\end{equation}
and the left-handed (${}_aD^{\beta}_{x}$) and the right-handed (${}_xD^{\beta}_{b}$)
space fractional derivatives in (\ref{eq:lin}) are the Riemann-Liouville fractional
derivatives of order $\beta$ \cite{AAKHMS,SSAKO} which are defined by
\begin{equation}
{}_aD^{\beta}_{x}u(x,t) = \frac{1}{\Gamma(2 - \beta)}\frac{\partial^2 }{\partial x^2}\int^{x}_a
\frac{u(\eta,t)d\eta}{(x -\eta)^{\beta - 1}}
\label{eq3x}
\end{equation}
and
\begin{equation}
{}_xD^{\beta}_{b}u(x,t) = \frac{1}{\Gamma(2 - \beta)}\frac{\partial^2 }{\partial x^2}
\int^{b}_x \frac{u(\eta,t)d\eta}{(\eta - x)^{\beta - 1}},
\label{eq3y}
\end{equation}
where $\Gamma(\cdot)$ denotes the Gamma function. Truly, when $\alpha = 1$ and $\beta = 2$, the
above equation reduces to the classical convection-diffusion equation (CDE).

The fractional CDE has been recently treated by a number of authors. It is presented as a
useful approach for the description of transport dynamics in complex systems which are
governed by anomalous diffusion and non-exponential relaxation patterns \cite{RMJKT,AISGMZ}.
The fractional CDE is also used in groundwater hydrology research to model the transport of
passive tracers carried by fluid flow in a porous medium \cite{DABSW}. Though analytic approaches,
such as the Fourier transform method, the Laplace transform methods, and the Mellin transform
method, have been proposed to seek closed-form solutions \cite{IPFDE,YPST,YZPF}, there are
very few FDEs whose analytical closed-form solutions are available. Therefore, the research on
numerical approximation and techniques for the solution of FDEs has attracted intensive interest.
Most early established numerical methods are developed for handling the space factional CDE or
the time fractional CDE. For space fractional CDE, many researchers exploited the conventional shifted
Gr\"{u}nwald discretization \cite{MMMCT} and the implicit Euler (or Crank-Nicolson) time-step discretization
for two-sided Riemman-Liouville fractional derivatives and the first order time derivative,
respectively. Then they constructed many numerical treatments for space fractional CDE, refer
to \cite{MMMCT,LSPC,LSWWZ,ESAS,ZDAXM,LSWWHW,KWHWA,ZDVSLB,FLPZKB,SMAAR} and
references therein for details. Later, Chen and Deng combined the second-order discretization
with the Crank-Nicolson temporal discretization for producing the novel numerical methods,
which archive the second order accuracy in both time and space for space fractional CDE \cite{MCWD,WDMCE}.
Even Chen \& Deng and Qu et al. separately designed the fast computational techniques, which
can also reduce the required algorithmic storage, for implementing the above mentioned second
order numerical method; see \cite{WDMCE,WQSLL} for details. Additionally, there are also some other
interesting numerical methods for space fractional CDE, refer, e.g., to
\cite{NJFLP,AHBDBA,AHBMAZ1,HHTMFLS,WYTWD2} for a discussion of these issues.

On the other hand, for time fractional CDE, many early implicit numerical methods are derived
by the combination of the $L1$ approximate formula \cite{YLCXF} for Caputo fractional derivative
and the first/second order spatial discretization. These numerical methods are unconditionally
convergent with the $\mathcal{O}(\tau^{2- \alpha} + h)$ or $\mathcal{O}(\tau^{2- \alpha} + h^2)$
accuracy, where $\tau$ and $h$ are time-step size and spatial grid size, respectively. In order
to improve the spatial accuracy, Cui \cite{MCAH,MCCES} and Mohebbi \& Abbaszadeh \cite{AMMAC}
proposed, respectively, two compact exponential methods and a compact finite difference method
for a time fractional convection-subdiffusion equation so that the spatial accuracy is improved
to the fourth-order. However, the methods and analyses in \cite{MCAH,AMMAC} are only for the
equations with constant coefficients. In particular, the discussions in \cite{AMMAC} are limited
to a special time fractional convection-subdiffusion equation where the diffusion and convection
coefficients are assumed to be one. In addition, some other related numerical methods have already
been proposed for handling the time fractional CDE, see e.g. \cite{SMAA,ZWSV,ZJFWC,SZXFYH,PZYTG,YMWAC} for
details.

On contrast, although the numerical methods for space (or time) fractional CDE are extensively
investigated in the past research, the work about numerically handling the TSFCDE is not too
much. Firstly, Zhang \cite{YZAFD,YZFDA}, Shao \& Ma \cite{YSWM}, Qin \& Zhang \cite{PQXZA},
and Liu et al. \cite{FLPZKB} have worked out a series of studies about constructing the implicit
difference scheme (IDS) for TSFCDE, however all these numerical scheme can archive the convergence
with first order accuracy in both space and time from the both theoretical and numerical
perspectives. Moreover, Liu et al. \cite{FLPZV,FLPZKB}, Zhao et al. \cite{ZZXQJ}, and Shen
\cite{SSFLVA} had considered to solve the more general form of TSFCDE, in which the first-order
space derivative is replaced by the two-sided Riemman-Liouville fractional derivative of order
$\nu\in (0, 1)$. Again, their numerical methods cannot enjoy the convergence with second order
accuracy in both space and time. In addition, some other efficient approaches are also developed
for dealing with TSFCDE numerically. Moreover, most of these numerical methods have no complete
theoretical analysis for both convergence and stability; e.g., refer to
\cite{MPMRE,YCYWY,SIMDS,AHBMAZ,AHBMAZ1,HHTMFL,WJYLA,JWYCBL} for details.

Traditional methods for solving FDEs tend to generate full coefficient matrices, which incur
computational cost of $\mathcal{O}(N^3)$ and storage of $O(N^2)$ with $N$ being the number of
grid points \cite{KWHWA,WQSLL}.
To optimize the computational complexity, a shifted Gr\"{u}nwald discretization scheme with
the property of unconditional stability was proposed by Meerschaet and Tadjeran \cite{MMMCT}
to approximate the space fractional CDE. Later, Wang and Wang \cite{KWHWA} discovered that the
linear system generated by this discretization has a special Toeplitz-like coefficient matrix,
or, more precisely, this coefficient matrix can be expressed as a sum of diagonal-multiply-Toeplitz
matrices. This implies that the storage requirement is $\mathcal{O}(N)$ instead of $\mathcal{O}(N^2)$,
and the complexity of the matrix-vector multiplication only requires $\mathcal{O}(N \log N)$ operations
by the fast Fourier transform (FFT) \cite{MNIM}. Upon using this advantage, Wang and Wang
proposed the CGNR method having computational cost of $\mathcal{O}(N \log^2 N)$ to solve the linear
system and numerical experiments show that the CGNR method is fast when the diffusion coefficients
are very small, i.e., the discretized systems are well-conditioned \cite{KWHWA}. However, the
discretized systems become ill-conditioned when the diffusion coefficients are not small. In this
case, the CGNR method converges slowly. To remedy this shortcoming, Zhao et al. have extended the preconditioned
technique, which is introduced by Lin et al. in the context of space fractional diffusion equation \cite{FRLSWY},
for handling the Toeplitz-like linear systems coming from numerical discretization of TSFCDE \cite{ZZXQJ}. Their
results related to the promising acceleration of the convergence of the iterative methods, while solving
(\ref{eq:lin}).

In this paper, we firstly derive an implicit difference scheme for solving (\ref{eq:lin}) and then
verify that the proposed scheme can archive the stability and convergence with second order accuracy
in both space and time (i.e. $\mathcal{O}(\tau^2 + h^2)$) from the both theoretical and numerical
perspectives. As far as we know, the scheme is the first one (without extrapolation techniques) who can have the convergence with $
\mathcal{O}(\tau^2 + h^2)$. For the time marching of the scheme, plenty of linear systems, which
possess different Toeplitz coefficient matrices, are required to be solved. Those linear systems
can be solved efficiently one after one by using Krylov subspace method with suitable circulant
preconditioners \cite{MNIM,SLLHWS}, then it can reduce the computational cost and memory deeply.
Especially for TSFCDE with constant coefficients, we turn to represent the inverse of the Toeplitz
coefficient matrix as a sum of products of Toeplitz triangular matrices \cite{MNIM,IGAS}, so that the solution
of each linear system for time marching can be obtained directly via several fast Fourier transforms
(FFTs). To obtain the explicit representation of Toeplitz matrix inversion, only two specific linear
systems with the same Toeplitz coefficient matrix are needed to be solved, which can be done by using
the preconditioned Krylov subspace methods \cite{YSMHS} with complexity $O(N\log N)$.

An outline of this paper is as follows. In Section \ref{sec2}, we establish
an implicit difference scheme for (\ref{eq:lin}) and prove that this scheme
is unconditionally stable and convergent with $O(\tau^2 + h^2)$ order accuracy.
In Section \ref{sec3}, we investigate that the resulting linear systems have
the nonsymmetric Toeplitz matrices, then we design the fast solution techniques
based on preconditioned Krylov subspace methods to solve (\ref{eq:lin}) by
exploiting the Toeplitz matrix property of the implicit difference scheme.
Finally, we present numerical experiments to illustrate the effectiveness of
our numerical approaches in Section \ref{sec4} and provide concluding remarks
in Section \ref{sec5}.

\section{Implicit difference scheme}
\label{sec2}
\quad\
In this section, we present an implicit difference method for discretizing the
TSFCDE defined by (\ref{eq:lin}). Unlike former numerical approaches with
the first order accuracy in both time and space \cite{FLPZV,ZZXQJ,YZFDA,YZAFD,YSWM},
we exploit henceforth two-sided fractional derivatives to approximate the Riemann-Liouville
derivatives in (\ref{eq3x}) and (\ref{eq3y}). We can show that, by two-sided fractional
derivatives, this proposed method is also unconditionally stable and convergent under
second order accuracy in time and space.

\subsection{Numerical discretization of the TSFCDE}
\quad\
In order to derive the proposed scheme, we first introduce some notations. In the rectangle
$\bar{Q}_T = \{(x,t): a \leq x \leq b, 0 \leq t \leq T \}$ we introduce the mesh $\varpi_{h\tau}
= \varpi_h \times \varpi_{\tau}$, where $\varpi_h = \{x_i = ih,~i = 0, 1,\ldots,N;~hN = b - a\}$,
and $\varpi_{\tau} = \{t_j = j\tau,~j = 0,1,\ldots,M;~\tau = T/M\}$. Besides, ${\bm
v} = \{v_i\mid 0\leq i \leq N\}$ be any grid function. Then, the following lemma introduced in
\cite{4A2015} gives a description on the time discretization.

\begin{lemma}
Suppose $0 < \alpha < 1$, $\sigma = 1 - \frac{\alpha}{2}, u(t) \in \mathcal{C}^{3}[0, T ]$,
and $t_{j + \sigma} = (j + \sigma)\tau$. Then
\begin{equation*}
\partial^{\alpha}_{0,t_{j + \sigma}}u(t) - \Delta^{\alpha}_{0,t_{j + \sigma}}u(t) =
\mathcal{O}(\tau^{3 - \alpha}),
\end{equation*}
where $\Delta^{\alpha}_{0,t_{j + \sigma}}u(t) = \frac{\tau^{-\alpha}}{\Gamma(2 - \alpha)}
\sum^{j}_{s = 0} c^{(\alpha,\sigma,j)}_{j - s}[u(t_{s+1}) - u(t_s)]$, and $c^{(\alpha,\sigma,
0)}_0 = a^{(\alpha, \sigma)}_0$ for $j = 0$,
\begin{equation*}
c^{(\alpha,\sigma,j)}_m =
\begin{cases}
a^{(\alpha, \sigma)}_0 + b^{(\alpha, \sigma)}_1, & m = 0,\\
a^{(\alpha, \sigma)}_m + b^{(\alpha, \sigma)}_{m + 1} - b^{(\alpha, \sigma)}_m, & 1\leq m \leq j-1,\\
a^{(\alpha, \sigma)}_j - b^{(\alpha, \sigma)}_j, & m = j,\\
\end{cases}
\end{equation*}
for $j \geq 1$, in which $a^{(\alpha, \sigma)}_0 = \sigma^{1 - \alpha}$, $a^{(\alpha, \sigma)}_{\ell}
= (\ell + \sigma)^{1 - \alpha} - (\ell - 1 + \sigma )^{1 - \alpha},$ for $\ell \geq 1$;
and $b^{(\alpha, \sigma)}_{\ell} = \frac{1}{2 - \alpha}[(\ell + \sigma)^{2 - \alpha} - (\ell - 1 +
\sigma)^{2 - \alpha}] - \frac{1}{2} [(\ell + \sigma)^{1 - \alpha} + (\ell - 1 + \sigma)^{1 - \alpha}]$.
\label{lem1}
\end{lemma}

Denote $\mathfrak{L}^{n + \beta}(\mathbb{R}) = \{v|v \in L_1(\mathbb{R})\ \mathrm{and}\
\int^{+\infty}_{-\infty} (1 + |k|)^{n + \beta}|\hat{v}(k)|dk < \infty\}$, where $\hat{v}(k)
= \int^{+\infty}_{-\infty} e^{\iota kx} v(x)dx$ is the Fourier transformation of $v(x)$,
and we use $\iota = \sqrt{-1}$ to denote the imaginary unit. For the discretization on space,
we introduce the following Lemma:
\begin{lemma} (\cite{ZPHZZS,SVPLX}) Suppose that $v \in \mathfrak{L}^{n + \beta}(\mathbb{R})$, and let
\begin{equation*}
\delta^{\beta}_{x,+} v(x) = \frac{1}{h^{\beta}}\sum^{[\frac{x - a}{h}]}_{k = 0}\omega^{(\beta)}_k
v(x - (k - 1)h),
\end{equation*}
\begin{equation*}
\delta^{\beta}_{x,-} v(x) = \frac{1}{h^{\beta}}\sum^{[\frac{b - x}{h}]}_{k = 0}\omega^{(\beta)}_k
v(x + (k - 1)h),
\end{equation*}
then for a fixed $h$, we have
\begin{equation*}
{}_aD^{\beta}_{x}v(x) = \delta^{\beta}_{x,+} v(x) + \mathcal{O}(h^2),
\end{equation*}
\begin{equation*}
{}_xD^{\beta}_{b}v(x) = \delta^{\beta}_{x,-} v(x) + \mathcal{O}(h^2),
\end{equation*}
where
\begin{equation*}
\begin{cases}
\omega^{(\beta)}_0 = \lambda_1g^{(\beta)}_0,\quad\ \omega^{(\beta)}_1 = \lambda_1g^{(\beta)}_1 +
\lambda_0g^{(\beta)}_0,\\
\omega^{(\beta)}_k = \lambda_1g^{(\beta)}_k + \lambda_0g^{(\beta)}_{k - 1} + \lambda_{-1}
g^{(\beta)}_{k-2},& k \geq 2,
\end{cases}
\end{equation*}
with
\begin{equation*}
\lambda_1 = \frac{\beta^2 + 3\beta + 2}{12},\quad \lambda_0 = \frac{4 - \beta^2}{6},\quad
\lambda_{-1} = \frac{\beta^2 - 3\beta + 2 }{12},\quad and\ \ g^{(\beta)}_k = (-1)^k
\binom{\beta}{k}.
\end{equation*}
\label{lem2}
\end{lemma}

Let $u(x,t) \in \mathcal{C}^{4,3}_{x,t}$ be a solution of the problem (\ref{eq:lin}). Let us consider
Eq. (\ref{eq:lin}) for $(x,t) = (x_i,t_{j + \sigma})\in \bar{Q}_T,~i = 1, 2,\ldots, N - 1,~
j = 0, 1,\ldots, M-1,~\sigma = 1 - \frac{\alpha}{2}$:
\begin{equation*}
\begin{split}
\partial^{\alpha}_{0,t_{j + \sigma}}u(x,t) = &~\gamma(t_{j + \sigma})\Big(\frac{\partial u(x,t)}{\partial x}
\Big)_{(x_i, t_{j + \sigma})} + d_{+}(t_{j + \sigma})\Big({}_aD^{\beta}_{x}u(x,t)\Big)_{(x_i, t_{j +
\sigma})} + \\
& d_{-}(t_{j + \sigma})\Big({}_xD^{\beta}_{b}u(x,t)\Big)_{(x_i, t_{j + \sigma})} + f(x_i, t_{j + \sigma}).
\end{split}
\end{equation*}
For simplicity, let us introduce some notations
\begin{equation*}
u^{(\sigma)}_i = \sigma u^{j + 1}_i + (1 - \sigma) u^{j}_i,\quad \gamma^{(j + \sigma)} = \gamma(t_{j + \sigma}),
\quad D^{(j + \sigma)}_{\pm} = d_{\pm}(t_{j + \sigma}),\quad f^{j + \sigma}_i = f(x_i, t_{j + \sigma})
\end{equation*}
and
\begin{equation*}
\delta^{\beta}_{h}u^{(\sigma)}_i = \gamma^{(j + \sigma)} \frac{u^{(\sigma)}_{i + 1} - u^{(\sigma)}_{i - 1}}{2h}
+ \frac{D^{(j + \sigma)}_+}{h^{\beta}}\sum^{i + 1}_{k = 0}\omega^{(\beta)}_{k} u^{(\sigma)}_{i - k + 1} +
\frac{D^{(j + \sigma)}_-}{h^{\beta}}\sum^{N - i + 1}_{k = 0}\omega^{(\beta)}_{k} u^{(\sigma)}_{i + k - 1}.
\end{equation*}
Then with regard to Lemma \ref{lem1} we derive the implicit difference scheme with the approximation
order $\mathcal{O}(h^2 + \tau^2)$:
\begin{equation}
\begin{cases}
\Delta^{\alpha}_{0,t_{j + \sigma}}u_i = \delta^{\beta}_{h}u^{(\sigma)}_i + f^{j + \sigma}_i, &
1\leq i \leq N - 1,\quad 0\leq j \leq M - 1,\\
u^{0}_i = \phi(x_i), & 1\leq i \leq N - 1, \\
u^{j}_0 = u^{j}_N = 0, & 0 \leq j \leq M.
\end{cases}
\label{ids1}
\end{equation}
It is interesting to note that for $\alpha \rightarrow1$ we obtain the Crank-Nicolson difference scheme.

\subsection{Analysis of the implicit difference scheme}
\quad\
In this section, we analyze the stability and convergence for the implicit difference scheme (\ref{ids1}).
Let
\begin{equation*}
V_h = \{{\bm v}\mid {\bm v} = \{v_i\}\ {\rm is}\ {\rm a}\ {\rm grid}\ {\rm function}\
{\rm on}\ \varpi_h\ {\rm and}\ v_i = 0\ {\rm if}\ i = 0,N\}
\end{equation*}
For $\forall {\bm u},~{\bm v}\in V_h$, we define the discrete inner product and the
corresponding discrete $L_2$ norm as follows,
\begin{equation*}
({\bm u},{\bm v}) = h\sum^{N-1}_{i = 1}u_i v_i,~~{\rm and}~~\|{\bm u}\|
= \sqrt{({\bm u}, {\bm u})}.
\end{equation*}
Now, some lemmas are provided for proving the stability and convergence of implicit
difference scheme (\ref{ids1}).
\begin{lemma} (\cite{ZPHZZS,MMMCT,FRLSWY})
Let $1 < \alpha < 2$ and $g^{(\beta)}_k$ be defined in Lemma \ref{lem2}. Then we have
\begin{equation*}
\begin{cases}
g^{(\beta)}_0 = 1,\quad g^{(\beta)}_1 = -\beta,\quad g^{(\beta)}_2 > g^{(\beta)}_3 > \cdots > 0,\\
\sum^{\infty}_{k = 0}g^{(\beta)}_k = 0,\quad \sum^{N}_{k = 0}g^{(\beta)}_k < 0,\quad\ N > 1,\\
g^{(\beta)}_k = \mathcal{O}(k^{-(\beta + 1)}),\quad\ g^{(\beta)}_k = \Big(1 - \frac{\beta + 1}{k}
\Big)g^{(\beta)}_{k -1},\quad k = 1,2,\ldots.
\end{cases}
\end{equation*}
\label{lem3}
\end{lemma}

\begin{lemma} (\cite{ZPHZZS,SVPLX})
Let $1 < \alpha < 2$ and $g^{(\beta)}_k$ be defined in Lemma \ref{lem2}. Then we have
\begin{equation*}
\begin{cases}
\omega^{(\beta)}_0 = 1,\quad \omega^{(\beta)}_1 < 0,\quad \omega^{(\beta)}_k > 0,\quad k \geq 3,\\
\sum^{\infty}_{k = 0}\omega^{(\beta)}_k = 0,\quad \sum^{N}_{k = 0}\omega^{(\beta)}_k < 0,\quad\ N > 1,\\
\omega^{(\beta)}_0  + \omega^{(\beta)}_2 \geq 0.
\end{cases}
\end{equation*}
\label{lem4}
\end{lemma}

\begin{lemma} (\cite{ZPHZZS,SVPLX})
For $1 <\beta < 2$, and any ${\bm v} \in V_h$, it holds that
\begin{equation*}
(\delta^{\beta}_{x,+}{\bm v},{\bm v}) = (\delta^{\beta}_{x,-}{\bm v},{\bm v})
\leq \Big(\frac{1}{h^{\beta}}\sum^{N - 1}_{k = 0}\omega^{(\beta)}_k\Big)\|{\bm v}
\|^2.
\end{equation*}
\label{lem5}
\end{lemma}
\begin{lemma}
For $1 < \beta < 2$, $N \geq 5$, and any ${\bm v}\in V_h$, there exists a positive
constant $c_1$, such that
\begin{equation*}
(-\delta^{\beta}_{x,+}{\bm v},{\bm v}) = (-\delta^{\beta}_{x,-}{\bm v},{\bm v})
> c_1\ln 2\|{\bm v}\|^2.
\end{equation*}
\label{lem6}
\end{lemma}
\noindent\textbf{Proof}. Since
\begin{equation*}
\sum^{2N + 2 }_{k = N}\omega^{(\beta)}_k = \sum^{2N}_{k = N}g^{(\beta)}_k
+ (\lambda_1 + \lambda_0)g^{(\beta)}_{2N + 1} + \lambda_1g^{(\beta)}_{2N
+ 2} + \zeta(\beta),
\end{equation*}
where
\begin{equation*}
\begin{split}
\zeta(\beta) = (\lambda_0 + \lambda_{-1})g^{(\beta)}_{N - 1} + \lambda_{-1}
g^{(\beta)}_{N - 2} & = \Big[(\lambda_0 + \lambda_{-1})\frac{N - 2 - \beta}
{N - 1} + \lambda_{-1}\Big]g^{(\beta)}_{N - 2} \\
& = \frac{(12 - 6\beta)N + \beta^3 + 4\beta^2 - \beta - 22}{12(N - 1)}
g^{(\beta)}_{N - 2},\\
& \triangleq \frac{\vartheta(\beta)}{12(N - 1)}
g^{(\beta)}_{N - 2}
\end{split}
\end{equation*}
then $\zeta(2) = 0$, $\vartheta(2) = 0$ and $\vartheta'(\beta) = -6N + 3\beta^2
+ 8\beta - 1\leq 27 - 6N$, which implies $\zeta(\beta)$ is a decreasing function
for $\beta\in [1, 2]$, if $N \geq 5$ and $\vartheta'(\beta)< 0$. Hence $\zeta(\beta)
> 0$ when $N \geq 5$.

Then, by Lemma \ref{lem3}, there exist positive constants $\tilde{c}_1$ and $c_1$,
such that
\begin{equation}
\begin{split}
\frac{1}{h^\beta}\sum^{\infty}_{k = N}\omega^{(\beta)}_k  & > \frac{1}{h^\beta}\sum^{
2N + 2}_{k = N}\omega^{(\beta)}_k > \frac{1}{h^\beta}\sum^{2N}_{k = N}g^{(\beta)}_k
\geq \tilde{c}_1\sum^{2N}_{k = N}k^{-(\beta + 1)}N^{\beta} \\
& > \tilde{c}_1\sum^{2N}_{k = N}
k^{-(\beta + 1)}\Big(\frac{k}{2}\Big)^{\beta}\\
& = \frac{\tilde{c}_1}{2^{\beta}}\sum^{2N}_{k = N}\frac{1}{k} \geq c_1\int^{2N}_N
\frac{1}{x}dx = c_1\ln 2,\quad\ N \geq 5.
\end{split}
\label{eq2.3}
\end{equation}
Using Lemmas \ref{lem4} and \ref{lem5}, we then obtain
\begin{equation*}
(-\delta^{\beta}_{x,+}{\bm v},{\bm v}) = (-\delta^{\beta}_{x,-}{\bm v},{\bm v})
> \Big(\frac{1}{h^{\beta}}\sum^{\infty}_{k = N}\omega^{(\beta)}_k\Big)\|{\bm v}\|^2
> c_1\ln 2\|{\bm v}\|^2,
\end{equation*}
which proves the lemma. \hfill $\Box$

Based on the above lemmas, we can obtain the following theorem, which is essential for
analyzing the stability of the proposed implicit difference scheme.

\begin{theorem}
For any ${\bm v}\in V_h$, it holds that
\begin{equation*}
(\delta^{\beta}_h {\bm v}, {\bm v}) \leq -c \ln2\|{\bm v}\|^2,
\end{equation*}
where $c$ is a positive constant independent of the spatial step size $h$.
\label{thm2.1}
\end{theorem}
\noindent\textbf{Proof}. The concrete expression of $(\delta^{\beta}_h {\bm v}, {\bm v})$ can
be written by
\begin{equation}
(\delta^{\beta}_h {\bm v}, {\bm v}) = \gamma^{(j + \sigma)}h\sum^{N - 1}_{i = 1} \frac{v_{i + 1} -
v_{i - 1}}{2h}v_i + D^{(j + \sigma)}_+ (\delta^{\beta}_{x,+}{\bm v},{\bm v}) + D^{(j + \sigma)}_{-}
(\delta^{\beta}_{x,-}{\bm v},{\bm v}).
\label{eq2.22}
\end{equation}
It notes that $v^0 = v^N = 0$, we have
\begin{equation}
\gamma^{(j + 1)}h\sum^{N - 1}_{i = 1} \frac{v_{i + 1} - v_{i - 1}}{2h}v_i = 0.
\label{eq2.23}
\end{equation}

Moreover, according to Lemma \ref{lem6}, there exists a positive constant $c_1$ independent
of the spatial step size $h$, such that for any non-vanishing vector ${\bm u}\in V_h$, we obtain
\begin{equation}
D^{(j + \sigma)}_+ (\delta^{\beta}_{x,+}{\bm v},{\bm v}) + D^{(j + \sigma)}_{-}
(\delta^{\beta}_{x,-}{\bm v},{\bm v}) \leq - c_1\ln2 \Big(D^{(j + \sigma)}_+ + D^{(j + \sigma)}_-\Big)
\|{\bm v}\|^2
\label{eq2.24}
\end{equation}
Let $c = c_1\Big(D^{(j + \sigma)}_+ + D^{(j + \sigma)}_{-}\Big)$. Inserting (\ref{eq2.23}) and (\ref{eq2.24})
into (\ref{eq2.22}), Theorem \ref{thm2.1} holds. \hfill $\Box$

\begin{lemma} (\cite{4A2015,SVPLX})
Let $V_{\tau} = \{{\bm u}| {\bm u} = (u^0, u^1,\ldots, u^M)\}$ For any ${\bm u}\in
V_{\tau}$; one has the following inequality
\begin{equation*}
[\sigma u^{j + 1} + (1 - \sigma) u^{j}] \Delta^{\alpha}_{0, t_{j + \sigma}}{\bm u}
\geq \frac{1}{2}\Delta^{\alpha}_{0, t_{j + \sigma}}({\bm u})^2.
\end{equation*}
\end{lemma}

Now we can conclude the stability and convergence of the implicit difference scheme (\ref{ids1}). For simplicity
of presentation, in our proof, we denote $a^{j+1}_s = \frac{c^{(\alpha,\sigma,j)}_{j - s}}
{\tau^{\alpha} \Gamma(2 - \alpha)}$. Then $\Delta^{\alpha}_{0, t_{j + \sigma}}u = \sum^{j
}_{s = 0} (u^{s + 1} - u^s)a^{j + 1}_s$.

\begin{theorem}
Denote $\|f^{j + \sigma}\|^2 = h\sum^{N - 1}_{i = 1}f^2(x_i,t_{j + \sigma})$. Then the
implicit difference scheme (\ref{ids1}) is unconditionally stable and the following a priori
estimate holds:
\begin{equation}
\|u^{j + 1}\|^2 \leq \|u^0 \|^2 + \frac{T^{\alpha}\Gamma(1 - \alpha)}{2c\ln 2}\|f^{j
 + \sigma}\|^2,\quad\ 0 \leq j \leq M - 1,
\label{add3x}
\end{equation}
where $u^{j + 1} = (u^{j+1}_1,u^{j+1}_2,\ldots,u^{j+1}_{N - 1})^T$.
\label{them2.2}
\end{theorem}
\textbf{Proof}. To make an inner product of (\ref{ids1}) with $u^{(\sigma)}$, we have
\begin{equation}
(\Delta^{\alpha}_{0,t_{j + \sigma}}u, u^{(\sigma)}) = (\delta^{\beta}_{h}u^{(\sigma)},
u^{(\sigma)}) + (f^{j + \sigma},u^{(\sigma)}).
\label{eq2.1xd}
\end{equation}
It follows from Theorem \ref{thm2.1} and Lemma \ref{lem6} that
\begin{eqnarray}
(\delta^{\beta}_h u^{(\sigma)}, u^{(\sigma)}) \leq -c \ln2\|u^{(\sigma)}\|^2, \label{eq2.30}\\
(\Delta^{\alpha}_{0,t_{j + \sigma}}u, u^{(\sigma)}) \geq \frac{1}{2} \Delta^{\alpha}_{0,t_{j +
\sigma}}(\| u\|^2). \label{eq2.31}
\end{eqnarray}
Inserting (\ref{eq2.30})-(\ref{eq2.31}) into (\ref{eq2.1xd}) and using the Cauchy-Schwarz
inequality, we obtain
\begin{equation}
\begin{split}
\frac{1}{2} \Delta^{\alpha}_{0,t_{j + \sigma}}(\| u\|^2) &\leq -c \ln2\|u^{(\sigma)}\|^2
+ (f^{j + \sigma},u^{(\sigma)})\\
&\leq -c \ln2\|u^{(\sigma)}\|^2 + c \ln2\|u^{(\sigma)}\|^2 + \frac{1}{8c\ln 2}\|f^{j +
\sigma}\|^2   \\
&\leq  \frac{1}{8c\ln 2}\|f^{j + \sigma}\|^2.
\end{split}
\label{eq2.32}
\end{equation}

Next, it holds that
\begin{equation*}
a^{j + 1}_j \|u^{j+1}\|^2 \leq \sum^{j}_{s = 1}(a^{j + 1}_s - a^{j + 1}_{s-1})\|u^s\|^2
+ a^{j + 1}_0 \|u^0 \|^2 + \frac{1}{4c\ln 2}\|f^{j + \sigma}\|^2.
\end{equation*}
Then, to notice that $a^{j + 1}_0 > \frac{1}{2T^{\alpha}\Gamma(1 - \alpha)}$ (cf. \cite{4A2015}),
we obtain
\begin{equation}
a^{j + 1}_j \|u^{j+1}\|^2 \leq \sum^{j}_{s = 1}(a^{j + 1}_s - a^{j + 1}_{s-1})\|u^s\|^2
+ a^{j + 1}_0 \Big(\|u^0 \|^2 + \frac{T^{\alpha}\Gamma(1 - \alpha)}{2c\ln 2}\|f^{j + \sigma}\|^2 \Big).
\label{add2x}
\end{equation}

Suppose $h < 1$ and denote
\begin{equation*}
\check{\mathcal{P}} \triangleq \|u^0 \|^2 + \frac{T^{\alpha}\Gamma(1 - \alpha)}{2c\ln 2}\|f^{j + \sigma}\|^2
\end{equation*}
Then, Eq. (\ref{add2x}) can be rewritten by
\begin{equation}
a^{j + 1}_j \|u^{j+1}\|^2 \leq \sum^{j}_{s = 1}(a^{j + 1}_s - a^{j + 1}_{s-1})\|u^s\|^2
+ a^{j + 1}_0 \check{\mathcal{P}}.
\label{add2y}
\end{equation}
Next, we prove that the estimate relation (\ref{add3x}) is valid for $j = 0,1\ldots,M - 1$ by mathematical
induction. Obviously, It follows from (\ref{add2y}) that (\ref{add3x}) holds for $j = 0$. Let us assume that
the inequality (\ref{add3x}) takes place for all $0 \leq j \leq k~(0 \leq k \leq M - 1)$, that is
\begin{equation*}
\|u^{j}\| \leq \check{\mathcal{P}},\qquad j = 0,1,\ldots,k.
\end{equation*}
For (\ref{add2y}) at $j = k$, one has
\begin{equation*}
\begin{split}
a^{k + 1}_k \|u^{k + 1}\|^2 & \leq \sum^{k}_{s = 1}(a^{k + 1}_s - a^{k + 1}_{s - 1})\|u^s\|^2
+ a^{k + 1}_0 \check{\mathcal{P}}  \\
& \leq \sum^{k}_{s = 1}(a^{k + 1}_s - a^{k + 1}_{s - 1})\check{\mathcal{P}}
+ a^{k + 1}_0 \check{\mathcal{P}}  = a^{k + 1}_k \check{\mathcal{P}}.
\end{split}
\end{equation*}
The proof of Theorem \ref{them2.2} is fully completed. \hfill $\Box$
\begin{theorem}
Suppose that $u(x,t)$ is the solution of (\ref{eq:lin}) and $\{u^{j}_{i}\mid x_i\in
\varpi_h,~~0\leq j\leq M\}$, is the solution of the implicit difference scheme (\ref{ids1}).
Denote
\begin{equation*}
E^{j}_i = u(x_i,t_j) - u^{j}_i,\quad x_i\in \varpi_h, \quad 0 \leq j\leq M.
\end{equation*}
Then there exists a positive constant $\tilde{c}$ such that
\begin{equation*}
\|E^j\| \leq \tilde{c} (\tau^2 + h^2),\quad\ 0 \leq j \leq M.
\end{equation*}
\end{theorem}
\textbf{Proof}. It can easily obtain that $E^j$ satisfies the following error equation
\begin{equation*}
\begin{cases}
\Delta^{\alpha}_{0,t_{j + \sigma}}E_i = \delta^{\beta}_{h}E^{(\sigma)}_i = R^{j + \sigma}_i, &
1\leq i \leq N - 1,\quad 0\leq j \leq M - 1,\\
E^{0}_i = 0, & 1\leq i \leq N - 1, \\
E^{j}_0 = E^{j}_N = 0, & 0 \leq j \leq M.
\end{cases}
\end{equation*}
where $R^{j + \sigma}_i = \mathcal{O}(\tau^2 + h^{2})$. By exploiting Theorem \ref{them2.2}, we
get
\begin{equation*}
\|E^{j + 1}\|^2 \leq \frac{T^{\alpha}\Gamma(1 - \alpha)}{2c\ln 2}\|R^{j + \sigma}\|^2
\leq \tilde{c} (\tau^2 + h^2),\quad\ 0 \leq j \leq M - 1,
\end{equation*}
which proves the theorem.  \hfill $\Box$

\section{Fast solution techniques based on preconditioned iterative solvers}
\label{sec3}
\quad\
In the section, we contribute to establish the efficient methods for solving a group of linear systems
with Toeplitz coefficient matrices, which are arisen from the matrix form of the implicit difference
scheme (\ref{ids1}). First of all, we derive the essential matrix form of the implicit difference
scheme (\ref{ids1}). Using notations in Section \ref{sec2}, the coefficient matrix of (\ref{ids1})
corresponding to each time level $j$ can be written as the following form,
\begin{equation}
\begin{cases}
\Big[\eta_j I - \sigma \Big(\frac{\gamma^{(\sigma)}}{2h}Q + \frac{D^{(\sigma)}_{+}}
{h^{\beta}}W_{\beta} + \frac{D^{(\sigma)}_{-}}{h^{\beta}}W^{T}_{\beta} \Big)\Big]{\bm u}^{1} = \Big[\eta_j I +
(1 - \sigma) \Big(\frac{\gamma^{(\sigma)}}{2h}Q +  \frac{D^{(\sigma)}_{+}} {h^{\beta}}W_{\beta}\\
+ \frac{D^{(\sigma)}_{-}}{h^{\beta}}W^{T}_{\beta} \Big)\Big]{\bm u}^{0} +
{\bm f}^{\sigma},\qquad\qquad\quad\ \ j = 0,\\
\Big[\eta_j I - \sigma \Big(\frac{\gamma^{(j + \sigma)}}{2h}Q + \frac{D^{(j +
\sigma)}_{+}}{h^{\beta}}W_{\beta} + \frac{D^{(j + \sigma)}_{-}}{h^{\beta}}W^{T}_{\beta}
\Big)\Big]{\bm u}^{j + 1} = \Big[\eta_j  I + (1 - \sigma) \Big(
\frac{\gamma^{(j + \sigma)}}{2h}Q~+ \\
\frac{D^{(j + \sigma)}_{+}}{h^{\beta}}W_{\beta} + \frac{D^{(j + \sigma)}_{-}}{h^{\beta}}W^{T}_{\beta} \Big)\Big]{\bm u}^{j}
- \frac{\tau^{-\alpha}}{\Gamma(2 - \alpha)} \sum\limits^{j - 1}_{s = 0} c^{(\alpha,\sigma,j)}_{j - s}({\bm u}^{s + 1} -
{\bm u}^s) + {\bm f}^{j + \sigma},\\
\qquad\qquad\qquad\qquad\qquad\qquad\qquad\ \ j = 1,2,\ldots,M-1,
\end{cases}
\label{eq3.1}
\end{equation}
where we have the coefficient
\begin{equation*}
\eta_j = \begin{cases}
\frac{c^{(\alpha,\sigma,0)}_0}{\tau^{\alpha}\Gamma(2 - \alpha)} = \frac{a^{(\alpha,\sigma)}_0}
{\tau^{\alpha}\Gamma(2 - \alpha)},& j =0\\
\frac{c^{(\alpha,\sigma,j)}_0}{\tau^{\alpha}\Gamma(2 - \alpha)} = \frac{a^{(\alpha,\sigma)}_0 + b^{(\alpha,\sigma)}_1}
{\tau^{\alpha}\Gamma(2 - \alpha)},& j = 1,2,\ldots,M - 1,
\end{cases}
\end{equation*}
then $Q$ and $W$ are two $(N - 1)\times (N - 1)$ real matrices of the following forms
\begin{equation}
Q = \begin{bmatrix}
0      & 1     & 0      & \cdots & 0      \\
-1     & 0     & 1      & \ddots & \vdots \\
0      & -1    & \ddots & \ddots & 0      \\
\vdots &\ddots & \ddots & \ddots & 1      \\
0      &\cdots & 0      & -1     & 0
\end{bmatrix},\quad\ \ \
W_{\beta} = \begin{bmatrix}
\omega^{(\beta)}_1       & \omega^{(\beta)}_0      & 0      & \cdots & 0      \\
\omega^{(\beta)}_2       & \omega^{(\beta)}_1      &\omega^{(\beta)}_0 & \ddots & \vdots \\
\vdots                   & \omega^{(\beta)}_2      &\omega^{(\beta)}_1& \ddots & 0      \\
\omega^{(\beta)}_{N - 2} & \cdots &\ddots          & \ddots &\omega^{(\beta)}_0    \\
\omega^{(\beta)}_{N - 1} &\omega^{(\beta)}_{N - 2} & 0      &\omega^{(\beta)}_2   & \omega^{(\beta)}_1
\end{bmatrix}.
\label{eq3.2}
\end{equation}
It is obvious that $W_{\beta}$ is a Toeplitz matrix (see \cite{MNIM,KWHWA,SVPLX}). Therefore, it can be stored
with $N + 1$ entries \cite{MNIM}.
\subsection{Resulting problems from the discretized scheme}
\quad\
According to (\ref{eq3.1}) and (\ref{eq3.2}), it indicates that the implicit difference scheme
(\ref{ids1}) requires to solve a nonsymmetric
Toeplitz linear system in each time level $j$, more precisely, there is a sequence of nonsymmetric
Toeplitz linear systems
\begin{equation}
A^{(j + \sigma)}{\bm u}^{(j + 1)} = B^{(j + \sigma)}{\bm u}^{(j)} + \delta{\bm u}^{(j)} + {\bm f}^{(j + \sigma)}
\label{eq3.3}
\end{equation}
where $\delta{\bm u}^{(j)}$ is used to denote the calculation of $\frac{\tau^{-\alpha}}{\Gamma(2 - \alpha)}
\sum\limits^{j - 1}_{s = 0} c^{(\alpha,\sigma,j)}_{j - s}({\bm u}^{s + 1} - {\bm u}^s)$ and
\begin{eqnarray*}
A^{(j + \sigma)} = \eta_j I - \sigma\Big(\frac{\gamma^{(j + \sigma)}}{2h}Q + \frac{D^{(j + \sigma)}_{+}}
{h^{\beta}}W_{\beta} + \frac{D^{(j + \sigma)}_{-}}{h^{\beta}}W^{T}_{\beta}\Big),\\
B^{(j + \sigma)} = \eta_j I + (1 - \sigma)\Big(\frac{\gamma^{(j + \sigma)}}{2h}Q + \frac{D^{(j + \sigma)}_{+}}
{h^{\beta}}W_{\beta} + \frac{D^{(j + \sigma)}_{-}}{h^{\beta}}W^{T}_{\beta}\Big),
\end{eqnarray*}
$j = 0,1,\ldots,M - 1$ and $A^{(j + 1)}$ varies with $j$; ${\bm f}^{(j + \sigma)}\in \mathbb{R}^{N - 1}$ also varies with $j$.
Here it should highlight that the sequence of linear systems (\ref{eq3.3}) corresponds to the time-stepping
scheme (\ref{ids1}), which is inherently sequential, hence the sequence of linear systems (\ref{eq3.3})
is extremely difficult to parallelize over time.

On the other hand, it is remarkable that if the coefficients $\gamma(t) = \gamma$ and $d_{\pm}(t) =
d_{\pm}$, then the coefficient matrices
\begin{equation}
A^{(j + \sigma)}= \begin{cases}
A^{(\sigma)},& j = 0,\\
A,& j = 1,2,\ldots,M - 1.
\end{cases}
\label{eqgux}
\end{equation}
Moreover, it is highlighted that the
coefficient $\eta_j$ is a real constant, which does not vary with $j = 1,2,\ldots,M - 1$. In other words, $A^{(j +
\sigma)} = A = \eta_j I - \sigma \Big(\frac{\gamma}{2h}Q + \frac{d_{+}}{h^{\beta}}W_{\beta} + \frac{d_{-}}{h^{
\beta}}W^{T}_{\beta}\Big)$ is independent of $j = 1,\ldots,M - 1$, i.e. the coefficient matrix of (\ref{eq3.3}) is
unchanged in each time level ($j = 1,2,\ldots,M$) of the implicit difference scheme. In this case, if we still solve linear systems (\ref{eq3.3})
one by one, it should be not sensible. A natural idea for this case is to find
the inverse of the Toeplitz matrix $A$, i.e. ${\bm u} = A^{-1}{\bm f}^{(j + \sigma)}$. It means that we are interested in computing $A^{-1}$ once
for all. One option is to compute the inverse by some direct methods such as the LU decomposition
\cite[pp. 44-54]{GHGCFV}. However, Toeplitz matrix is often dense, and the computation of the inverse of a large
dense matrix is prohibitive, especially when the matrix is large. Fortunately, as $A$ is also a
Toeplitz matrix, we have the Gohberg-Semencul formula (GSF) \cite{MNIM,IGAS} for its inverse.
Indeed, the inverse of a Toeplitz matrix $A$ can be reconstructed from its first
and last columns. More precisely, denote by ${\bm e}_1, {\bm e}_{N-1}$ the first and the last column
of the $(N - 1)$-by-$(N - 1)$ identity matrix, and let ${\bm x} = [\xi_0, \xi_1,\ldots, \xi_{N - 2}]^T$
and ${\bm y} = [\eta_0, \eta_1,\ldots, \eta_{N - 2}]^T$ be the solutions of the following two Toeplitz
systems
\begin{equation}
A{\bm x} = {\bm e}_1\quad\ \mathrm{and}\quad\ A{\bm y} = {\bm e}_{N - 1}.
\label{eq3.4}
\end{equation}
If $\xi_0 \neq 0$, then the Gohberg-Semencul formula can be expressed as
\begin{equation}
\begin{split}
A^{-1} = &~ \frac{1}{\xi_0}\Biggl(
\begin{bmatrix}
\xi_0 & 0 &  \cdots & 0\\
\xi_1 & \xi_0 & \cdots &0\\
\vdots & \vdots &\ddots & \vdots \\
\xi_{N-2} & \xi_{N - 3} & \cdots &\xi_0\\
\end{bmatrix}\cdot\begin{bmatrix}
\eta_{N-2} & \eta_{N - 3} & \cdots &\eta_0\\
0 & \eta_{N-2} & \cdots &\eta_1\\
\vdots & \vdots &\ddots & \vdots \\
0 & 0 &\cdots &\eta_{N-2}\\
\end{bmatrix}\\
& - \begin{bmatrix}
0    & 0 &  \cdots & 0 \\
\eta_0 & \cdots & 0 & 0 \\
\vdots & \vdots &\ddots & \vdots \\
\eta_{N-3} & \cdots & \eta_0  & 0\\
\end{bmatrix}\cdot\begin{bmatrix}
0 & \xi_{N-2} & \cdots &\xi_1   \\
\vdots & \vdots &\ddots & \vdots \\
0 & 0  & \cdots  &\xi_{N-2} \\
0 & 0 &\cdots &0\\
\end{bmatrix}
\Biggr) = \frac{1}{\xi_0}(L_p R_p - L^{0}_p R^{0}_p),
\end{split}
\end{equation}
where $L_p,~L^{0}_p$ are both lower Toeplitz matrices, and $R_p,R^{0}_p$ are upper
Toeplitz matrices. Consequently, the Toeplitz matrix-vector multiplication $A^{-1}
{\bm f}^{(j + \sigma)}$ can be archived in several FFTs of length $N - 1$ \cite{MNIM}.
For convenience, the following fast algorithm can be applied to compute the product
of $A^{-1}$ and a vector ${\bm v}$.

\begin{algorithm}[!htpb]
\caption{Compute ${\bm z} = A^{-1}{\bm v}$}
\begin{algorithmic}[1]
\STATE Solve two linear systems in Eq. (\ref{eq3.4})
\STATE Compute ${\bm z}_1 = R^{0}_p{\bm v}$ and ${\bm z}_2 = R_p{\bm v}$ via FFTs
\STATE Compute ${\bm z}_3 = L^{0}_p{\bm z}_1$ and ${\bm z}_4 = L_p{\bm z}_2$ via FFTs
\STATE Compute ${\bm z} = \frac{1}{\xi_0}({\bm z}_4 - {\bm z}_3)$
\end{algorithmic}
\label{alg1x}
\end{algorithm}

In summary, we need to search some efficient solvers for the nonsymmetric resulted
Toeplitz linear systems, whether to solve (\ref{eq3.3}) or to solve (\ref{eq3.4}).
In next subsection, we will introduce how to build efficient preconditioned iterative
solvers for nonsymmetric Toeplitz linear systems.

\subsection{Fast implementation of IDS based on preconditioned iterative solvers}
\quad\
In this subsection, we discuss the detailed framework about implementing the proposed
implicit difference scheme (\ref{ids1}). For the sake of clarity, an algorithm for solving
the implicit difference scheme is given in Algorithm \ref{alg1}.
\begin{algorithm}[!htpb]
\caption{Practical implementation of IDS}
\begin{algorithmic}[1]
  \FOR{$j = 0,1,\ldots,M-1$,}
  \STATE Compute ${\bm g}^{(j + 1)} = B^{(j + \sigma)}{\bm u}^{(j)} + \delta{\bm u}^{(j)}
  + {\bm f}^{(j + \sigma)}$
  \STATE Solve $A^{(j + \sigma)}{\bm u}^{(j + 1)} = {\bm g}^{(j + 1)}$
  \ENDFOR
\end{algorithmic}
\label{alg1}
\end{algorithm}

From Algorithm \ref{alg1}, $M$ real linear systems are required to be solved. If a direct
solver (e.g., Gaussian elimination method \cite[pp. 33-44]{GHGCFV}) is applied, its complexity will be $\mathcal{O}(MN^3)$,
which is very expansive if $N$ is large. Note that in Steps 2-3 of Algorithm \ref{alg1},
the Toeplitz structure of the matrices $A^{(j + \sigma)}$ and $B^{(j + \sigma)}$ has not
been utilized when solving the linear system. Actually, the matrix-vector product $B^{
(j + \sigma)}{\bm u}^{(j)}$  in Step 2 can be evaluated by FFTs in $\mathcal{O}(N\log N)
$ operations, then fast Toeplitz iterative solvers have been studied intensively, see e.g.,
\cite{MNIM} and references therein. Recently, they have been applied for space fractional
diffusion equation \cite{FRLSWY,SVPLX,XMGTZH,ZZXQJ,KWHWA,SLLHWS}. With Krylov subspace method
with circulant preconditioners in \cite{SLLHWS}, the Toeplitz system from space fractional
diffusion equation can be solved with a fast convergence rate. In this case, it also remarked
that the algorithmic complexity of preconditioned Krylov subspace methods is only in
$\mathcal{O}(N\log N)$ arithmetic operations per iteration step.

Motivated by the above considerations, we first implement the matrix-vector multiplication
$B^{(j + \sigma)}{\bm u}^{(j)}$ in Step 1 of Algorithm \ref{alg1} by FFTs. Then solving
the Toeplitz linear systems $A^{(j + \sigma)}{\bm u}^{(j + 1)} = {\bm g}^{(j + 1)}$ in
Step 3 by Krylov subspace methods, e.g., the conjugate gradient squared (CGS) method
\cite[pp. 241-244]{YSMHS}, with the circulant preconditoner
\begin{equation}
P^{(j + \sigma)} = \eta_j I - \sigma\Big(\frac{\gamma^{(j + \sigma)}}{2h}s(Q) + \frac{
D^{(j + \sigma)}_{+}}{h^{\beta}}s(W_{\beta}) + \frac{D^{(j + \sigma)}_{-}}{h^{\beta}}
s(W^{T}_{\beta})\Big),
\label{preku}
\end{equation}
where $s(\cdot)$ means the well-known Strang circulant approximation of a given Toeplitz matrix
\cite{MNIM,SLLHWS}. High efficiency of Strang circulant preconditioner for space FDEs
has been verified in \cite{SLLHWS}. To make sure the preconditioner defined in (\ref{preku}) is well-defined, let
us illustrate that $P^{(j + \sigma)}$ are nonsingular,
\begin{lemma}
All eigenvalues of $s(W_{\beta})$ and $s(W^{T}_{\beta})$ fall inside the open disc
\begin{equation}
\{z\in\mathbb{C}: |z - \omega^{(\beta)}_1| < -\omega^{(\beta)}_1\}
\end{equation}
\end{lemma}
\textbf{Proof}. All the Gershgorin disc \cite[pp. 119-122]{YSMHS} of the circulant matrices $s(W_{\beta})$
and $s(W^{T}_{\beta})$ are centered at $-\omega^{(\beta)}_1 > 0$ with radius
\begin{equation}
r_N = \omega^{(\beta)}_0 + \sum^{\frac{\lfloor N\rfloor}{2}}_{k = 2}
\omega^{(\beta)}_k < \sum^{\infty}_{k = 0,k \neq 1}\omega^{(\beta)}_k
= -\omega^{(\beta)}_1.
\end{equation}
by the properties of the sequence $\{\omega^{(\beta)}_k\}$; refer to Lemma
\ref{lem3} and Lemma \ref{lem4}. \hfill $\Box$

\begin{remark}
It is worth to mention that:
\begin{itemize}
  \item[1.] The real parts of all eigenvalues of $s(W_{\beta})$ and $s(W^{T}_{
  \beta})$ are strictly negative for all $N$;
  \item[2.] The absolute values of all eigenvalues of $s(W_{\beta})$ and $s(W^{
  T}_{\beta})$ are bounded above by $2|\omega^{(\beta)}_1|$ for all $N$.

\end{itemize}
\label{remk1}
\end{remark}

As we know, a circulant matrix can be quickly diagonalized by the Fourier matrix $F$
\cite{MNIM,WQSLL,SLLHWS}. Then it follows that $s(Q) = F^{*}\Lambda_{q}F$,
$s(W_{\beta}) = F^{*}\Lambda_{\beta}F$, and $s(W^{T}_{\beta}) = F^{*}\bar{
\Lambda}_{\beta}F$, where $\bar{\Lambda}_{\beta}$ is the complex conjugate
of $\Lambda_{\beta}$. Decompose the circulant matrix $P^{(j + \sigma)} =
F^{*}\Lambda_p F$ with the diagonal matrix $\Lambda_p = \eta_j I - \sigma
\Big(\frac{\gamma^{(j + \sigma)}}{2h}\Lambda_{q} + \frac{D^{(j + \sigma)}_{
+}}{h^{\beta}}\Lambda_{\beta} + \frac{D^{(j + \sigma)}_{-}}{h^{\beta}}
\bar{\Lambda}_{\beta}\Big)$. Then $P^{(j + \sigma)}$ is invertible if all
diagonal entries of $\Lambda_p$ are nonzero. Moreover, we can obtain the
following conclusion about the invertibility of $P^{(j + \sigma)}$ in
(\ref{preku}).

\begin{theorem}
The circulant preconditioners $P^{(j + \sigma)}$ defined as in (\ref{preku})
are nonsingular.
\label{themx2}
\end{theorem}
\textbf{Proof}. First of all, we already know that $Q$ is a skew-symmetric
Toeplitz matrix, it also finds that $s(Q)$ is also a skew-symmetric circulant matrix, thus the real part
$\Lambda_{q}$ is equal to zero, i.e., $\mathrm{Re}([\Lambda_q]_{k,k}) = 0$. On the
other hand, by Part 1 of Remark \ref{remk1}, we have $\mathrm{Re}([\Lambda_{\beta
}]_{k,k}) < 0$. Noting that $\eta_j > 0$, $\sigma>0$, and $D^{(j + \sigma)}_{
\pm} \geq 0$, thus we obtain
\begin{equation*}
\mathrm{Re}([\Lambda_p]_{k,k}) = \eta_j - \sigma\Big(0 + \frac{D^{(j + \sigma)}_{
+}}{h^{\beta}}\mathrm{Re}([\Lambda_{\beta}]_{k,k}) + \frac{D^{(j + \sigma)}_{-}}{
h^{\beta}}\mathrm{Re}([\bar{\Lambda}_{\beta}]_{k,k})\Big) \neq 0,
\end{equation*}
for each $k = 1,2,\ldots, N - 1$. Therefore, $P^{(j + \sigma)}$ are invertible.
\hfill $\Box$

Here although we do not plan to theoretically investigate the eigenvalue distributions of preconditioned matrices $(P^{(j + \sigma)})^{-1}A^{(j
+ \sigma)}$, we still can give some figures to illustrate the clustering eigenvalue
distributions of several specified preconditioned matrices in next section. Furthermore,
our numerical experiments show that the iteration numbers always fluctuates between 6 and 10, so
we regard the complexity of solving the linear system in Step 3 as $\mathcal{O}(N \log N)$.
It also implies that the computational complexity for implementing the whole IDS is about
$\mathcal{O}(MN\log N)$ operations.

Beside, if we have the coefficients $\gamma(t) = \gamma$ and $d_{\pm}(t) = d_{\pm}$ in
Eq. (\ref{eq:lin}), the matrix $B^{(j + \sigma)}$ also has the similar form as $B^{(j
+ \sigma)}$ in Eq. (\ref{eqgux}), then we can simplify Algorithm \ref{alg1} as the
following Algorithm \ref{alg2}.

\begin{algorithm}[!htpb]
\caption{Practical implementation of IDS with constant coefficients}
\begin{algorithmic}[1]
  \FOR{$j = 0,1,\ldots,M-1$,}
  \IF{$j = 0$}
  \STATE Compute ${\bm g}^{(1)} = B^{(\sigma)}{\bm u}^{(0)} + {\bm f}^{(\sigma)}$
  \STATE Solve $A^{(\sigma)}{\bm u}^{(1)} = {\bm g}^{(1)}$
  \ELSE
  \STATE Compute ${\bm g}^{(j+1)} = B{\bm u}^{(j)} + \delta{\bm u}^{(j)} + {\bm f
  }^{(j + \sigma)}$
  \STATE Solve ${\bm u}^{(j+1)} = A^{-1}{\bm g}^{(j + 1)}$
  \ENDIF
  \ENDFOR
\end{algorithmic}
\label{alg2}
\end{algorithm}

Again, if a direct method is employed to solve the linear system in Step 4 and compute
matrix inverse $A^{-1}$ with the help of LU decomposition\footnote{For the given linear system $A{\bm x} = {\bm b}$, we
solve it by MATLAB as: $\mathtt{[L,U]=lu(A);~~x = U\backslash(L\backslash b);}$}, which can be reused
in Step 7 of Algorithm \ref{alg2}, but its complexity will be still $\mathcal{O}(MN^3)$,
which is very costly if $N$ is large. Observe that in Steps 3, 4, 6 and 7 of Algorithm
\ref{alg2}, the Toeplitz structure of those four matrices has not been utilized when solving
the linear system. Actually, two matrix-vector multiplications $B^{(\sigma)}{\bm u}^{(0)}$
and $B{\bm u}^{(j)}$ in Steps 3 and 6 can be evaluated by FFTs in $\mathcal{O}(N\log N)$
operations, then fast Toeplitz iterative solvers with suitable circulant preconditioners,
the Toeplitz system in Step 4 and (\ref{eq3.4}) can be solved with a fast convergence rate.
Here we can construct two circulant preconditioners defined as
\begin{eqnarray}
P^{(\sigma)} = \eta_0 I - \sigma\Big(\frac{\gamma}{2h}s(Q) + \frac{D_{+}}{h^{\beta}}s(W_{\beta})
+ \frac{D_{-}}{h^{\beta}}s(W^{T}_{\beta})\Big), \label{eqx2x}\\
P = \eta_j I - \sigma\Big(\frac{\gamma}{2h}s(Q) + \frac{D_{+}}{h^{\beta}}s(W_{\beta}) +
\frac{D_{-}}{h^{\beta}}s(W^{T}_{\beta})\Big), \label{eqx2xzz}
\end{eqnarray}
for the linear systems in Step 4 and (\ref{eq3.4}), respectively. Here the invertibility
of those above two circulant preconditioners introduced in (\ref{eqx2x})-(\ref{eqx2xzz}) can be similarly archived by using Theorem \ref{themx2}. Then
we can employ Algorithm
\ref{alg1x} to evaluate the matrix-vector multiplication $A^{-1}{\bm g}^{(j + 1)}$  in Step 7 of Algorithm
\ref{alg2}. Similarly, according to results in next section, we find that the iteration numbers
required by preconditioned Krylov subspace methods always fluctuate between 6 and 10. In this case, it also
is remarkable that the algorithmic complexity of preconditioned Krylov subspace methods is only $\mathcal{O
}(N\log N)$ at each iteration step. In conclusion, the total complexity for implementing
the IDS with constant coefficients is also in $\mathcal{O}(MN\log N)$ operations.

\section{Numerical results}
\label{sec4}
\quad\
In this section we first carry out some numerical experiments to
illustrate that our proposed IDS can indeed converge with the
second order accuracy in both space and time. At the same time,
some numerical examples are reported to show the effectiveness of
the fast solution techniques (i.e., Algorithms \ref{alg1x}-\ref{alg2})
designed in Section \ref{sec3}. For Krylov subspace method and
direct solver, we choose built-in functions for the preconditioned
CGS (PCGS) method,
LU factorization of MATLAB in Example 1 and MATLAB's
backslash in Example 2, respectively. For the CGS method with
circulant preconditioners, the stopping criterion of those methods
is $\|{\bm r}^{(k)}\|_2/\|{\bm r}^{(0)}\|_2 < 10^{-12}$, where
${\bm r}^{(k)}$ is the residual vector of the linear system after
$k$ iterations, and the initial guess is chosen as the zero vector.
All experiments were performed on a Windows 7 (64 bit) PC-Intel(R)
Core(TM) i5-3740 CPU 3.20 GHz, 8 GB of RAM using MATLAB 2014a with
machine epsilon $10^{-16}$ in double precision floating point arithmetic.
\vspace{2mm}

\noindent\textbf{Example 1}. In this example, we consider the equation (\ref{eq:lin}) on space
interval $[a, b]=[0, 1]$ and time interval $[0, T]=[0, 1]$ with diffusion coefficients $d_+(t)
= d_{+} = 0.8,~d_{-}(t) = d_{-} = 0.5$, convection coefficient $\gamma(t) = \gamma = -0.1$,
initial condition $u(x,0) = x^2(1 - x)^2$, and source term
\begin{equation*}
\begin{split}
f(x,t) = &~\frac{\Gamma(3 + \alpha)}{2}x^2(1 - x)^2 t^2 - (t^{2 + \alpha} + 1)\Big\{2\gamma
x(1-x)(1-2x) + \frac{\Gamma(3)}{\Gamma(3 - \beta)}[d_{+}x^{2 - \beta} \\
& + d_{-}(1 - x)^{2 - \beta}] - \frac{2\Gamma(4)}{\Gamma(4 - \beta)}[d_{+}x^{3 - \beta} + d_{-}(1 - x)^{3 - \beta
}] + \frac{\Gamma(5)}{\Gamma(5 - \beta)}[d_{+}x^{4 - \beta} \\
& + d_{-}(1 - x)^{4 - \beta}]\Big\}.
\end{split}
\end{equation*}
The exact solution of this example is $u(x, t) = (t^{2 + \alpha} + 1)x^2(1 - x)^2$. For the
finite difference discretization, the space step and time step are taken to be $h =
1/N$ and $\tau = h$, respectively. The errors ($E = U - u$) and convergence order (CO) in
the norms $\|\cdot\|_0$ and $\|\cdot\|_{\mathcal{C}(\bar{\omega}_{h\tau})}$, where $\|U
\|_{\mathcal{C}(\bar{\omega}_{h\tau})} = \max_{(x_i,t_j)\in \bar{\omega}_{h\tau}} |U|$,
are given in Tables \ref{tab1}-\ref{tab2}. Here these notations are used throughout this
section. Additionally, the performance of fast solution techniques presented in Section
\ref{sec3} for this example will be illustrated in Tables \ref{tab2y}-\ref{tab2z}. In the
following tables ``\texttt{Speed-up}" is defined as
\begin{equation*}
{\tt Speed\texttt{-}up} = \frac{\mathtt{Time1}}{\mathtt{Time2}}.
\end{equation*}
Obviously, when ${\tt Speed\texttt{-}up}>1$, it means that \texttt{Time2} needed by our
proposed method is more competitive than \texttt{Time1} required by Algorithm \ref{alg2}
with reusing LU decomposition in aspects of the CPU time elapsed.

\begin{table}[!htpb]
\caption{{\small $L_2$-norm and maximum norm error behavior versus grid size reduction when
$\tau = h$ and $\beta = 1.8$ in Example 1.}}
\centering
\begin{tabular}{crcccc}
\toprule
$\alpha$ & $h$  & $\max_{0\leq n \leq M} \|E^n\|_0$ & CO in $\|\cdot\|_0$ &$\|E\|_{\mathcal{C}(\bar{\omega}_{h\tau})}$
&CO in $\|\cdot\|_{\mathcal{C}(\bar{\omega}_{h\tau})}$\\
\hline
0.10 & 1/32  & 2.7954e-4 & --     & 4.0880e-4 & --     \\
     & 1/64  & 6.6775e-5 & 2.0657 & 9.8580e-5 & 2.0520 \\
     & 1/128 & 1.6010e-5 & 2.0603 & 2.3815e-5 & 2.0494 \\
     & 1/256 & 3.8514e-6 & 2.0556 & 5.7630e-6 & 2.0470 \\
0.50 & 1/32  & 2.6670e-4 & --     & 3.8874e-4 & --     \\
     & 1/64  & 6.3583e-5 & 2.0685 & 9.3590e-5 & 2.0544 \\
     & 1/128 & 1.5219e-5 & 2.0628 & 2.2573e-5 & 2.0518 \\
     & 1/256 & 3.6558e-6 & 2.0576 & 5.4539e-6 & 2.0492 \\
0.90 & 1/32  & 2.4972e-4 & --     & 3.6255e-4 & -- \\
     & 1/64  & 5.9441e-5 & 2.0708 & 8.7173e-5 & 2.0562 \\
     & 1/128 & 1.4206e-5 & 2.0650 & 2.0993e-5 & 2.0540 \\
     & 1/256 & 3.4078e-6 & 2.0596 & 5.0762e-6 & 2.0481 \\
0.99 & 1/32  & 2.5899e-4 & --     & 3.7959e-4 & --     \\
     & 1/64  & 6.2121e-5 & 2.0598 & 9.1923e-5 & 2.0460 \\
     & 1/128 & 1.4944e-5 & 2.0555 & 2.2275e-5 & 2.0450 \\
     & 1/256 & 3.6057e-6 & 2.0512 & 5.4042e-6 & 2.0433 \\
\bottomrule
\end{tabular}
\label{tab1}
\end{table}

\begin{table}[!htpb]
\caption{{\small $L_2$-norm and maximum norm error behavior versus $\tau$-grid size reduction
when $h = 1/1000$ and $\beta = 1.8$ in Example 1.}}
\centering
\begin{tabular}{cccccc}
\toprule
$\alpha$ & $\tau$ & $\max_{0\leq n \leq M} \|E^n\|_0$ & CO in $\|\cdot\|_0$ &$\|E\|_{\mathcal{C}(\bar{\omega}_{h\tau})}$
&CO in $\|\cdot\|_{\mathcal{C}(\bar{\omega}_{h\tau})}$\\
\hline
0.10 & 1/10 & 1.9209e-5 & --     & 3.0437e-5 & --     \\
     & 1/20 & 4.6741e-6 & 2.0390 & 7.4069e-6 & 2.0389 \\
     & 1/40 & 1.0134e-6 & 2.2054 & 1.6095e-6 & 2.2023 \\
0.50 & 1/10 & 1.2639e-4 & --     & 1.9985e-4 & --     \\
     & 1/20 & 3.1564e-5 & 2.0015 & 4.9914e-5 & 2.0014 \\
     & 1/40 & 7.7315e-6 & 2.0295 & 1.2232e-5 & 2.0288 \\
0.90 & 1/10 & 2.4927e-4 & --     & 3.9380e-4 & --     \\
     & 1/20 & 6.2151e-5 & 2.0039 & 9.8203e-5 & 2.0036 \\
     & 1/40 & 1.5356e-5 & 2.0170 & 2.4272e-5 & 2.0165 \\
0.99 & 1/10 & 2.7402e-4 & --     & 4.3269e-4 & --     \\
     & 1/20 & 6.8333e-5 & 2.0036 & 1.0791e-4 & 2.0035 \\
     & 1/40 & 1.6913e-5 & 2.0145 & 2.6714e-5 & 2.0142\\
\bottomrule
\end{tabular}
\label{tab2}
\end{table}
As seen from Table \ref{tab1}, it finds that as the number of the spatial subintervals and time steps is increased keeping $h = \tau$,
a reduction in the maximum error takes place, as expected and the convergence order of the approximate scheme is
$\mathcal{O}(h^2) = \mathcal{O}(\tau^2)$, where the convergence order is given by the formula: CO = $\log_{h_1/h_2}
\frac{\|E_1\|}{\|E_2\|}$ ($E_i$ is the error corresponding to $h_i$). On the other hand, Table \ref{tab2} illustrates
that if $h = 1/1000$, then as the number of time steps of our approximate scheme is increased, a reduction in the
maximum error takes place, as expected and the convergence order of time is $\mathcal{O}(\tau^2)$, where the convergence
order is given by the following formula: CO = $\log_{\tau_1/\tau_2}\frac{\|E_1\|}{\|E_2\|}$.

\begin{figure}[!htbp]
\centering
\includegraphics[width=3.1in,height=2.95in]{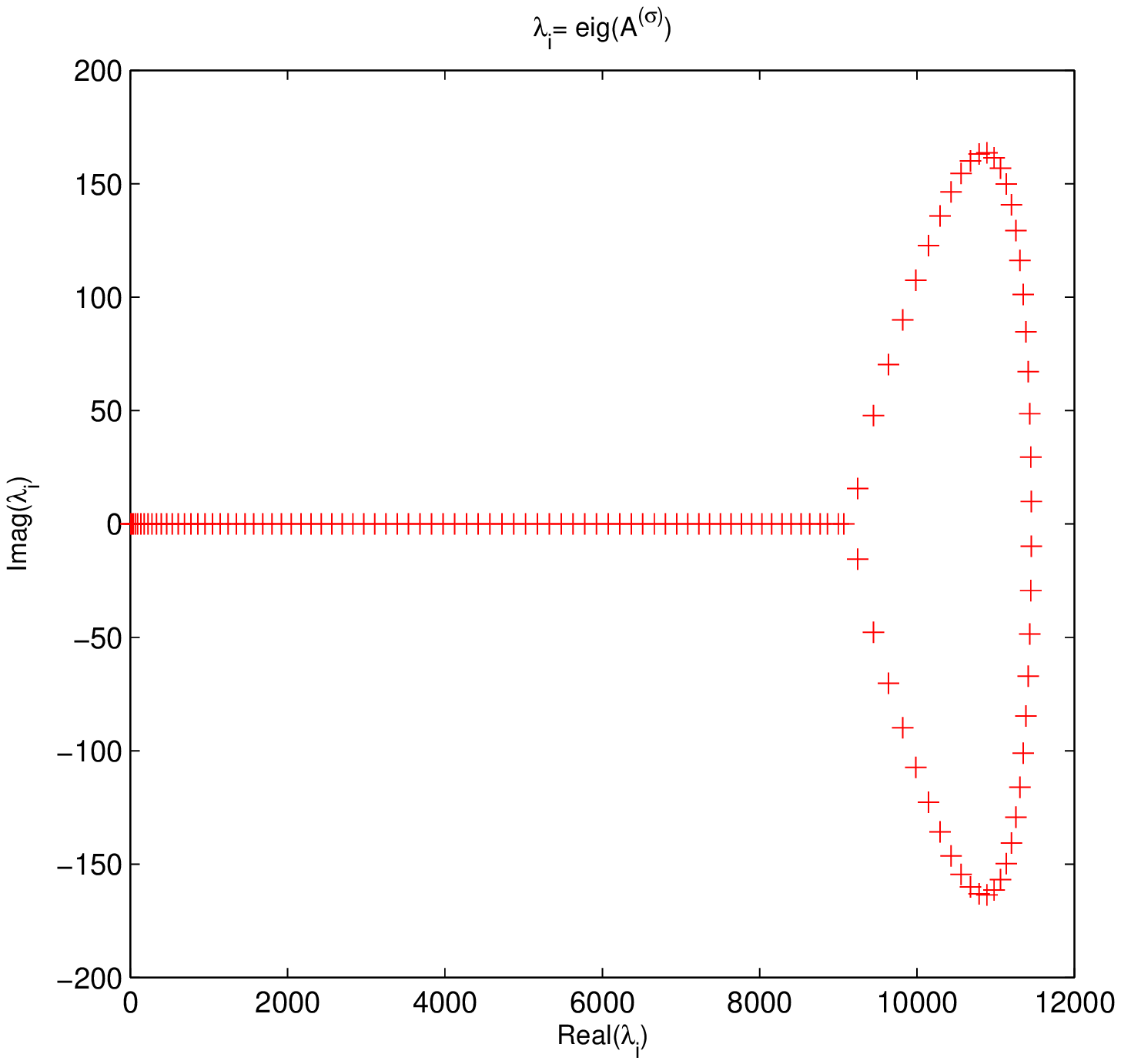}
\includegraphics[width=3.1in,height=2.95in]{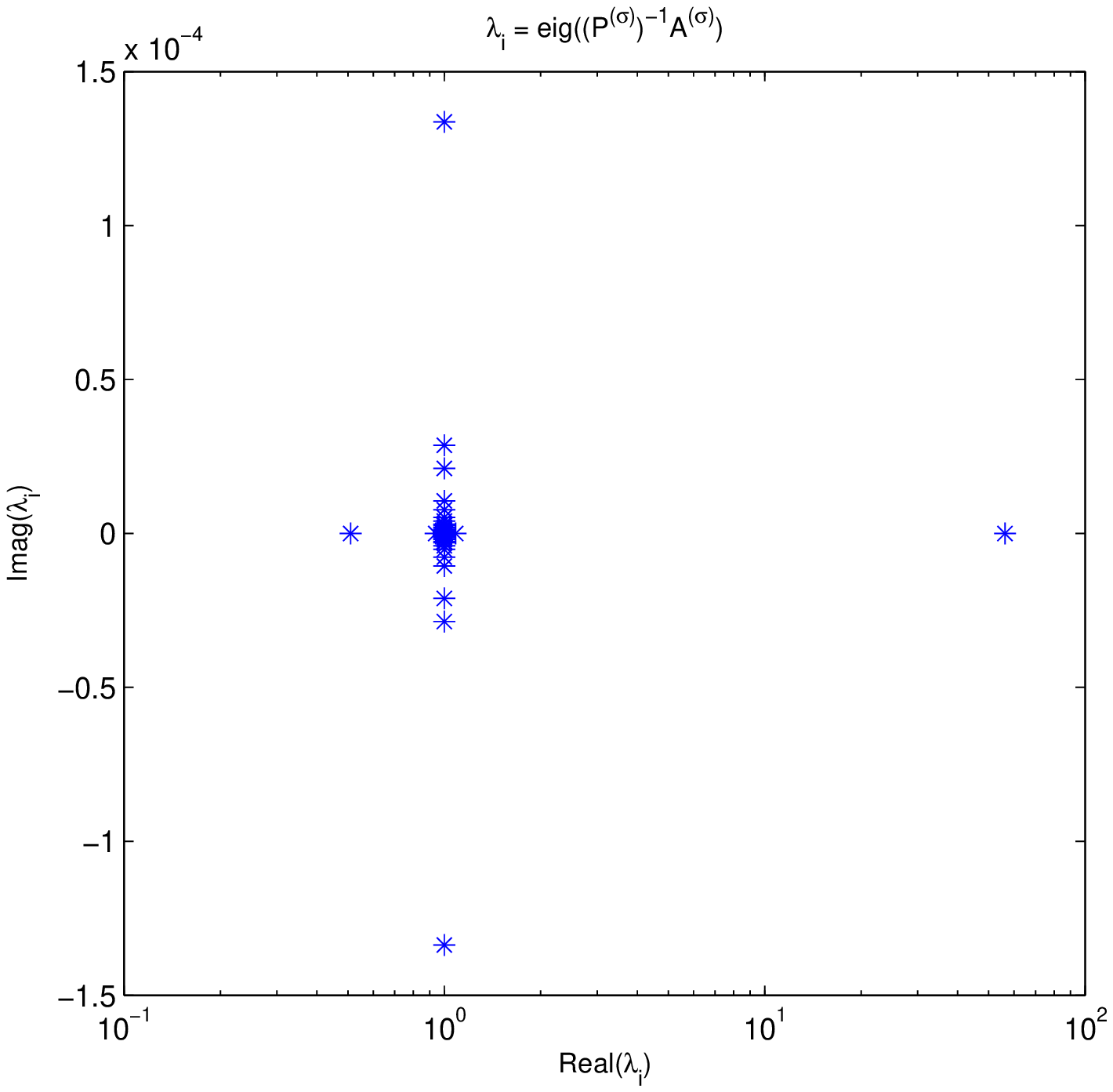}
\caption{{\small Spectrum of both original and preconditioned
matrices at the time level $j = 0$, respectively, when $N =
M = 128,~\alpha = 0.9$ and $\beta = 1.8$. Left: Original matrix; Right:
circulant preconditioned matrix.}}
\label{fig1}
\end{figure}

\begin{figure}[!htbp]
\centering
\includegraphics[width=3.1in,height=2.95in]{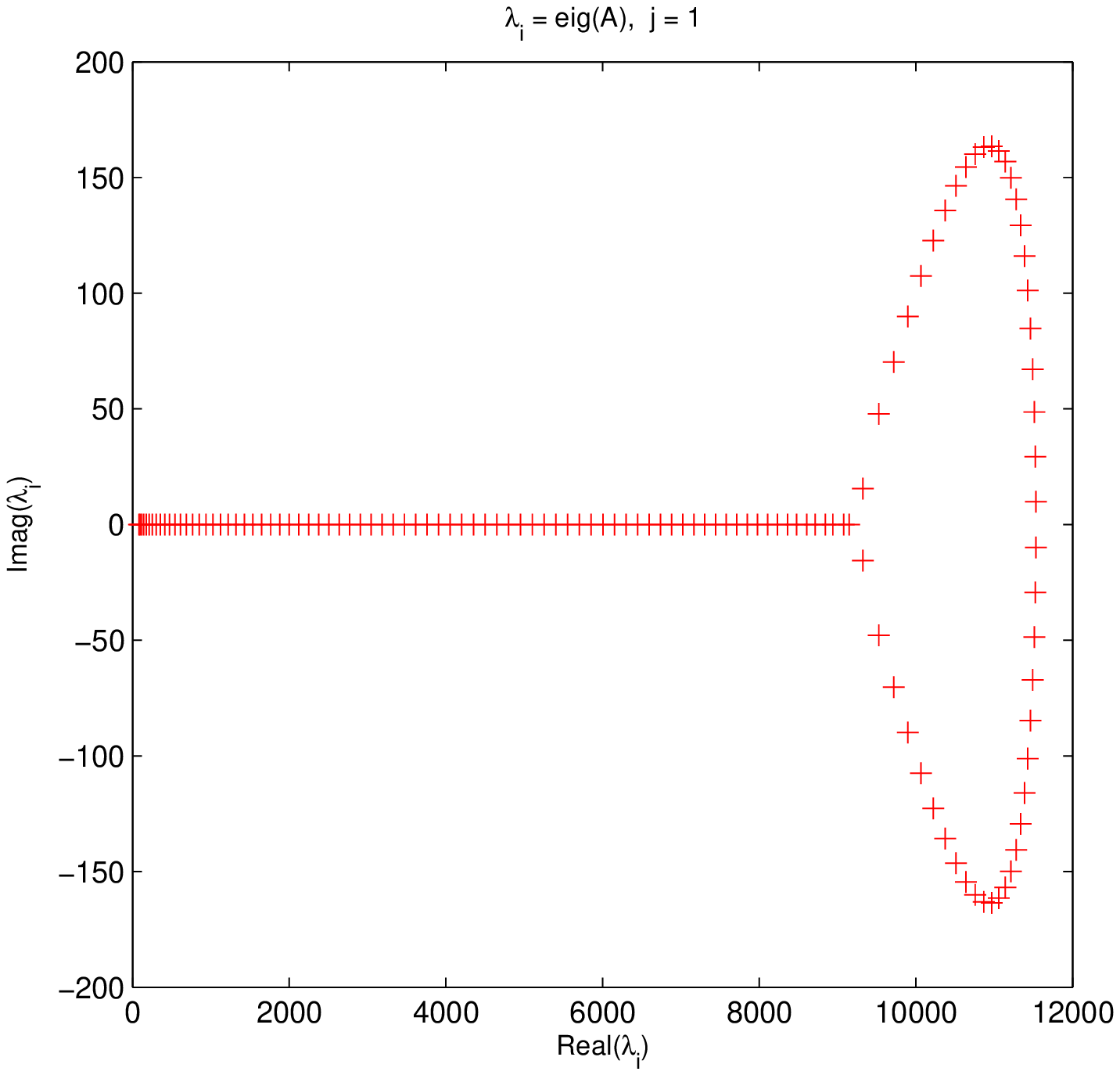}
\includegraphics[width=3.1in,height=2.95in]{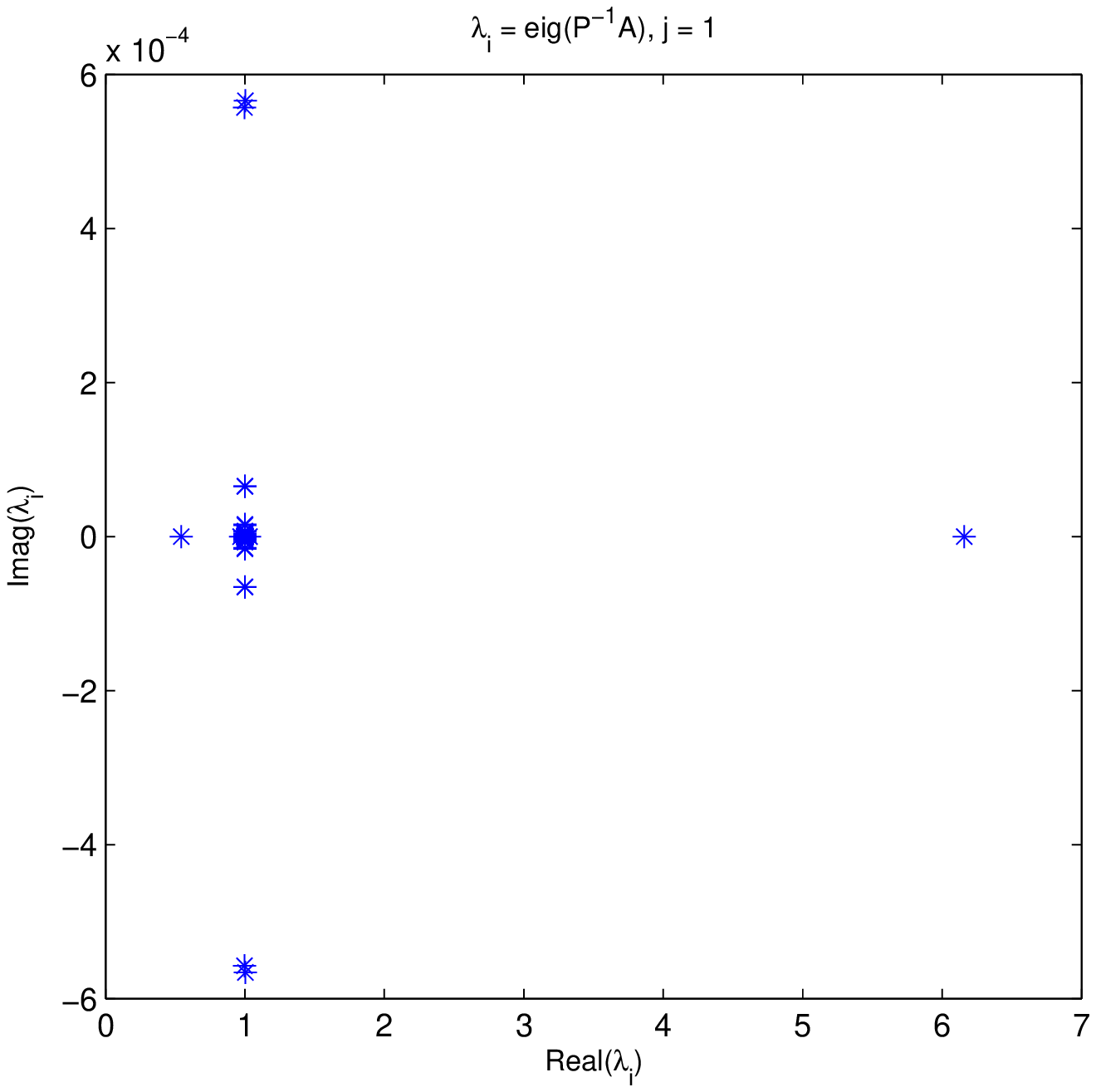}
\caption{{\small Spectrum of both original and preconditioned
matrices at the time level $j = 1$, respectively, when $N =
M = 128,~\alpha = 0.9$ and $\beta = 1.8$. Left: Original matrix;
Right: circulant preconditioned matrix.}}
\label{fig2}
\end{figure}

Firstly, some eigenvalue plots about both original and preconditioned matrices
are drawn in Figs. \ref{fig1}-\ref{fig2}. These two figures confirm that for
circulant preconditioning, the eigenvalues of preconditioned matrices are
clustered at 1, expect for few (about $6\sim 10$) outliers. The vast majority
of the eigenvalues are well separated away from 0. It may be interpreted as
that in our implementation the number of iterations required by preconditioned
Krylov subspace methods almost ranges from 6 to 10. We validate the effectiveness
and robustness of the designed circulant preconditioner from the perspective of
clustering spectrum distribution.

\begin{table}[t]\small\tabcolsep=5.3pt
\begin{center}
\caption{{\small CPU time in seconds for solving Example 1 with $\alpha = 0.9$ that
\texttt{Time1} is done by for Algorithm \ref{alg2} (LU decomposition) and \texttt{Time2}
is done by Algorithm \ref{alg2} with Algorithm \ref{alg1x}.}}
\begin{tabular}{lrrcrrcrrc}
\hline &\multicolumn{3}{c}{$\beta = 1.2$} &\multicolumn{3}{c}{$\beta = 1.5$} &\multicolumn{3}
{c}{$\beta = 1.8$}\\
[-2pt]\cmidrule(l{0.7em}r{0.7em}){2-4} \cmidrule(l{0.7em}r{0.6em}){5-7}\cmidrule(l{0.7em}r{0.7em}){8-10}\\[-11pt]
$h = \tau$ &\texttt{Time1} &\texttt{Time2} & $\texttt{Speed-up}$ &\texttt{Time1} &\texttt{Time2} &
$\texttt{Speed-up}$&\texttt{Time1} &\texttt{Time2} &$\texttt{Speed-up}$    \\
\hline
$2^{-5}$  &0.003  &0.009   &0.33   &0.003   &0.009  &0.33  &0.003  &0.009  &0.33  \\
$2^{-6}$  &0.011  &0.017   &0.65   &0.011   &0.017  &0.65  &0.011  &0.016  &0.69  \\
$2^{-7}$  &0.061  &0.059   &1.03   &0.062   &0.058  &1.05  &0.062  &0.058  &1.05  \\
$2^{-8}$  &0.344  &0.266   &1.29   &0.336   &0.267  &1.26  &0.335  &0.268  &1.25  \\
$2^{-9}$  &4.319  &1.902   &2.27   &4.317   &1.901  &2.27  &4.338  &1.899  &2.28  \\
$2^{-10}$ &41.683 &17.035  &2.43   &42.106  &17.556 &2.40  &41.808 &17.166 &2.43  \\
\hline
\end{tabular}
\label{tab2y}
\end{center}
\end{table}

\begin{table}[!htbp]\small\tabcolsep=5.3pt
\begin{center}
\caption{{\small CPU time in seconds for solving Example 1 with $\beta = 1.8$ that
\texttt{Time1} is done by for Algorithm \ref{alg2} (LU decomposition) and \texttt{Time2}
is done by Algorithm \ref{alg2} with Algorithm \ref{alg1x}.}}
\begin{tabular}{lrrcrrcrrc}
\hline &\multicolumn{3}{c}{$\alpha = 0.1$} &\multicolumn{3}{c}{$\alpha = 0.5$} &\multicolumn{3}
{c}{$\alpha = 0.9$}\\
[-2pt]\cmidrule(l{0.7em}r{0.7em}){2-4} \cmidrule(l{0.7em}r{0.6em}){5-7}\cmidrule(l{0.7em}r{0.7em}){8-10}\\[-11pt]
$h = \tau$ &\texttt{Time1} &\texttt{Time2} & $\texttt{Speed-up}$ &\texttt{Time1} &\texttt{Time2} &
$\texttt{Speed-up}$&\texttt{Time1} &\texttt{Time2} &$\texttt{Speed-up}$    \\
\hline
$2^{-5}$  &0.003  &0.010   &0.30   &0.003  &0.010  &0.30 &0.003  &0.010  &0.30 \\
$2^{-6}$  &0.011  &0.017   &0.65   &0.011  &0.017  &0.65 &0.011  &0.017  &0.65 \\
$2^{-7}$  &0.061  &0.059   &1.03   &0.060  &0.059  &1.02 &0.060  &0.058  &1.03 \\
$2^{-8}$  &0.343  &0.268   &1.28   &0.345  &0.268  &1.29 &0.340  &0.266  &1.28 \\
$2^{-9}$  &4.329  &1.903   &2.27   &4.349  &1.907  &2.28 &4.357  &1.900  &2.27 \\
$2^{-10}$ &41.639 &17.058  &2.44   &41.677 &17.097 &2.44 &41.662 &16.986 &2.45 \\
\hline
\end{tabular}
\label{tab2z}
\end{center}
\end{table}

In Tables \ref{tab2y}-\ref{tab2z}, it illustrates that the proposed fast direct
solver for different discretized problems takes much less CPU time elapsed as
$M$ and $N$ become large. When $M = N = 2^{10}$ and different discretized
parameters, the CPU time of Algorithm \ref{alg2} is about 17 seconds, the speedup is more than 2 times. Meanwhile,
although \texttt{Time1} required by Algorithm \ref{alg2} for small test
problems ($M = N = 32,64$) than \texttt{Time2} needed by Algorithm
\ref{alg2}, our proposed method is still more attractive in terms of
lower memory requirement. Compared to Algorithm \ref{alg2} with reusing
LU decomposition, it highlighted that in the whole implementation the proposed
solution technique does not require to store the full matrices (e.g. some
matrices $A^{(\sigma)},A$ and their LU factors) at all. In short, we can
conclude that our proposed IDS with fast implementation is still more
competitive than the IDS with reusing the conventional LU decomposition.
\vspace{2mm}

\noindent\textbf{Example 2}. In the last test, we investigate the equation (\ref{eq:lin})
on the space interval $[a,b] = [0,1]$ and the time interval $[0,T] = [0,1]$ with
diffusion coefficients $d_+(t) = 9\sin(t),~ d_{-}(t) = 4\sin(t)$, convection
coefficient $\gamma(t) = -t$, initial condition $u(x,0) = x^2(1 - x)^2
$, and source term
\begin{equation*}
\begin{split}
f(x,t) = & \frac{\Gamma(3 + \alpha)}{2}t^2 x^2(1-x)^2 - (t^{2 + \alpha} + 1)\Big\{-2tx(1-x)(1-2x)
+ \frac{\Gamma(3)\sin(t)}{\Gamma(3 - \beta)}[9x^{2 - \beta}\\
& + 4(1 - x)^{2 - \beta}] - \frac{2\Gamma(4)\sin(t)}{\Gamma(4 - \beta)}[9x^{3 - \beta} + 4(1 - x)^{3 - \beta
}] + \frac{\Gamma(5)\sin(t)}{\Gamma(5 - \beta)}[9x^{4 - \beta}~+\\
& 4(1 - x)^{4 - \beta
}]\Big\}.
\end{split}
\end{equation*}
The exact solution of this
example is defined as $u(x,t) = (t^{2 + \alpha})x^2(1-x)^2$. For the implicit finite
difference
discretization, the space step and time step are taken to be $h = 1/N$ and $\tau = h$,
respectively. The experiment results about the proposed IDS for Example 3 are reported in
Tables \ref{tab7}-\ref{tab8}.
Furthermore, the effectiveness of fast solution techniques presented in Section
\ref{sec3} for this example will be illustrated in Tables \ref{tab8y}-\ref{tab8z}.

\begin{table}[t]
\caption{{\small $L_2$-norm and maximum norm error behavior versus grid
size reduction when $\tau = h$ and $\beta = 1.3$ in Example 2.}}
\centering
\begin{tabular}{crcccc}
\toprule
$\alpha$ & $h$  & $\max_{0\leq n \leq M} \|E^n\|_0$ & CO in $\|\cdot\|_0$
&$\|E\|_{\mathcal{C}(\bar{\omega}_{h\tau})}$&CO in $\|\cdot\|_{\mathcal{C
}(\bar{\omega}_{h\tau})}$\\
\hline
0.10 & 1/32  & 3.1941e-4 & --     & 5.6886e-4 & --     \\
     & 1/64  & 7.6298e-5 & 2.0657 & 1.6055e-4 & 1.8250 \\
     & 1/128 & 1.8397e-5 & 2.0521 & 4.2694e-5 & 1.9110 \\
     & 1/256 & 4.4694e-5 & 2.0414 & 1.1036e-5 & 1.9519 \\
0.50 & 1/32  & 3.0866e-4 & --     & 5.6897e-4 & --     \\
     & 1/64  & 7.3673e-5 & 2.0668 & 1.6054e-4 & 1.8254 \\
     & 1/128 & 1.7757e-5 & 2.0527 & 4.2689e-5 & 1.9110 \\
     & 1/256 & 4.3137e-6 & 2.0414 & 1.1035e-5 & 1.9518 \\
0.90 & 1/32  & 2.9880e-4 & --     & 5.6951e-4 & --     \\
     & 1/64  & 7.1478e-5 & 2.0636 & 1.6058e-4 & 1.8264 \\
     & 1/128 & 1.7232e-5 & 2.0524 & 4.2691e-5 & 1.9113 \\
     & 1/256 & 4.1814e-6 & 2.0430 & 1.1034e-5 & 1.9519 \\
0.99 & 1/32  & 3.2304e-4 & --     & 5.7367e-4 & --     \\
     & 1/64  & 7.7278e-5 & 2.0638 & 1.6119e-4 & 1.8314 \\
     & 1/128 & 1.8633e-5 & 2.0522 & 4.2748e-5 & 1.9149 \\
     & 1/256 & 4.5227e-6 & 2.0426 & 1.1035e-5 & 1.9538 \\
\bottomrule
\end{tabular}
\label{tab7}
\end{table}

\begin{table}[!htbp]
\caption{{\small $L_2$-norm and maximum norm error behavior versus
$\tau$-grid size reduction when $h = 1/1200$ and $\beta = 1.3$ in
Example 2.}}
\centering
\begin{tabular}{cccccc}
\toprule
$\alpha$ & $\tau$ & $\max_{0\leq n \leq M} \|E^n\|_0$ & CO in $\|\cdot
\|_0$ &$\|E\|_{\mathcal{C}(\bar{\omega}_{h\tau})}$&CO in $\|\cdot\|_{
\mathcal{C}(\bar{\omega}_{h\tau})}$\\
\hline
0.10 & 1/10 & 2.0652e-5 & --     & 3.2538e-5 & --     \\
     & 1/20 & 5.0679e-6 & 2.0269 & 7.9617e-6 & 2.0310 \\
     & 1/40 & 1.1568e-6 & 2.1312 & 1.7951e-6 & 2.1490 \\
0.50 & 1/10 & 1.3380e-4 & --     & 2.1072e-4 & --     \\
     & 1/20 & 3.3465e-5 & 1.9994 & 5.2683e-5 & 1.9999 \\
     & 1/40 & 8.2623e-6 & 2.0180 & 1.2988e-5 & 2.0202 \\
0.90 & 1/10 & 2.6237e-4 & --     & 4.1300e-4 & --     \\
     & 1/20 & 6.5529e-5 & 2.0014 & 1.0313e-4 & 2.0016 \\
     & 1/40 & 1.6266e-5 & 2.0103 & 2.5583e-5 & 2.0112 \\
0.99 & 1/10 & 2.8754e-4 & --     & 4.5251e-4 & --     \\
     & 1/20 & 7.1771e-5 & 2.0023 & 1.1293e-4 & 2.0025 \\
     & 1/40 & 1.7825e-5 & 2.0095 & 2.8027e-5 & 2.0105 \\
\bottomrule
\end{tabular}
\label{tab8}
\end{table}

According to the numerical results illustrated in Table \ref{tab7}, it finds that as
the number of the spatial subintervals and time steps is increased keeping $h = \tau$,
a reduction in the maximum error takes place, as expected and the convergence order of
the approximate scheme is $\mathcal{O}(h^2) = \mathcal{O}(\tau^2)$, where the convergence
order is given by the formula: CO = $\log_{h_1/h_2}\frac{\|E_1\|}{\|E_2\|}$ ($E_i$ is
the error corresponding to $h_i$). On the other hand, Table \ref{tab8} illustrates that
if $h = 1/1000$, then as the number of time steps of our approximate scheme is increased,
a reduction in the maximum error takes place, as expected and the convergence order of
time is $\mathcal{O}(\tau^2)$, where the convergence order is given by the following
formula: CO = $\log_{\tau_1/\tau_2}\frac{\|E_1\|}{\|E_2\|}$.

\begin{figure}[!htbp]
\centering
\includegraphics[width=3.1in,height=2.95in]{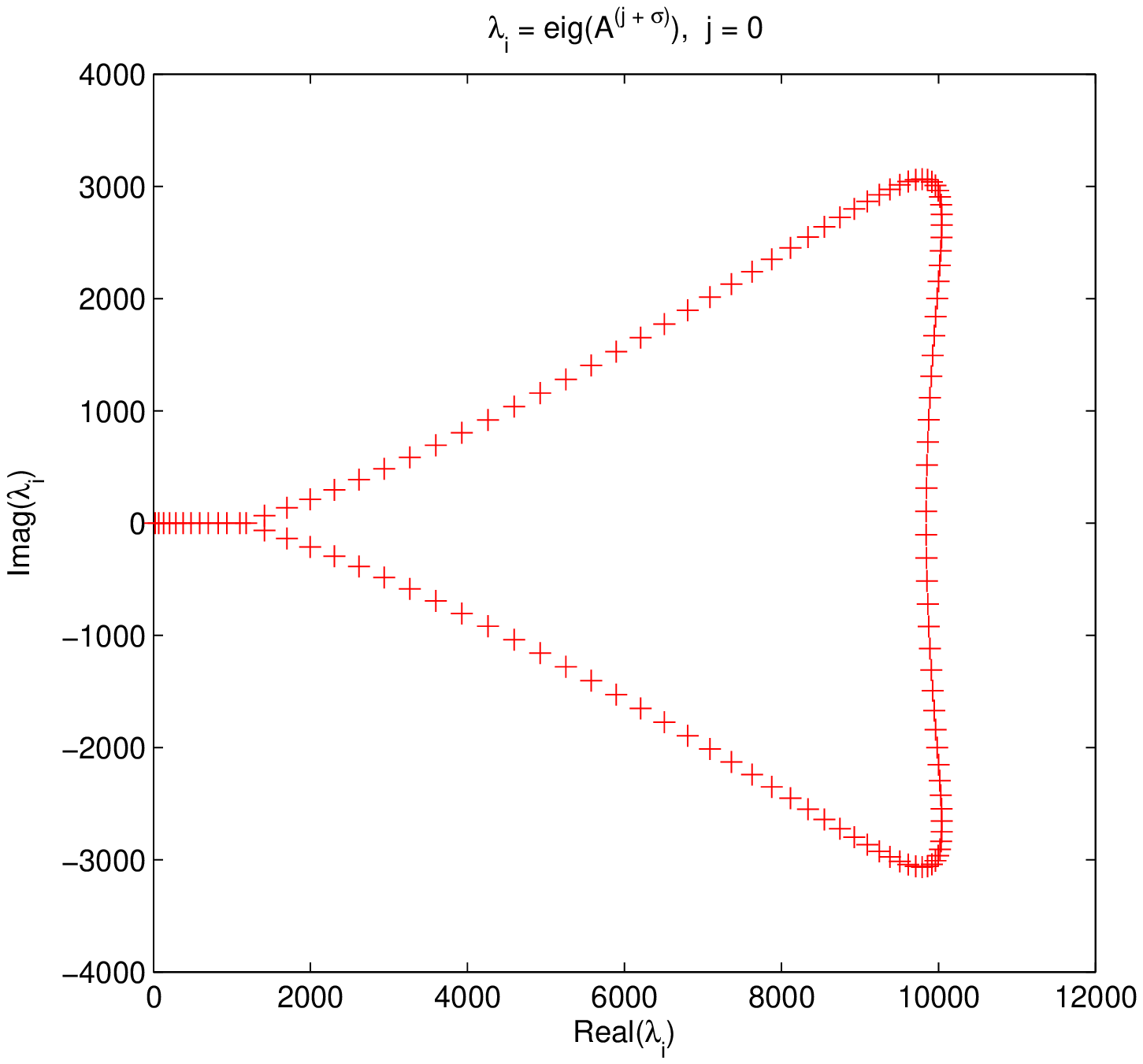}
\includegraphics[width=3.1in,height=2.95in]{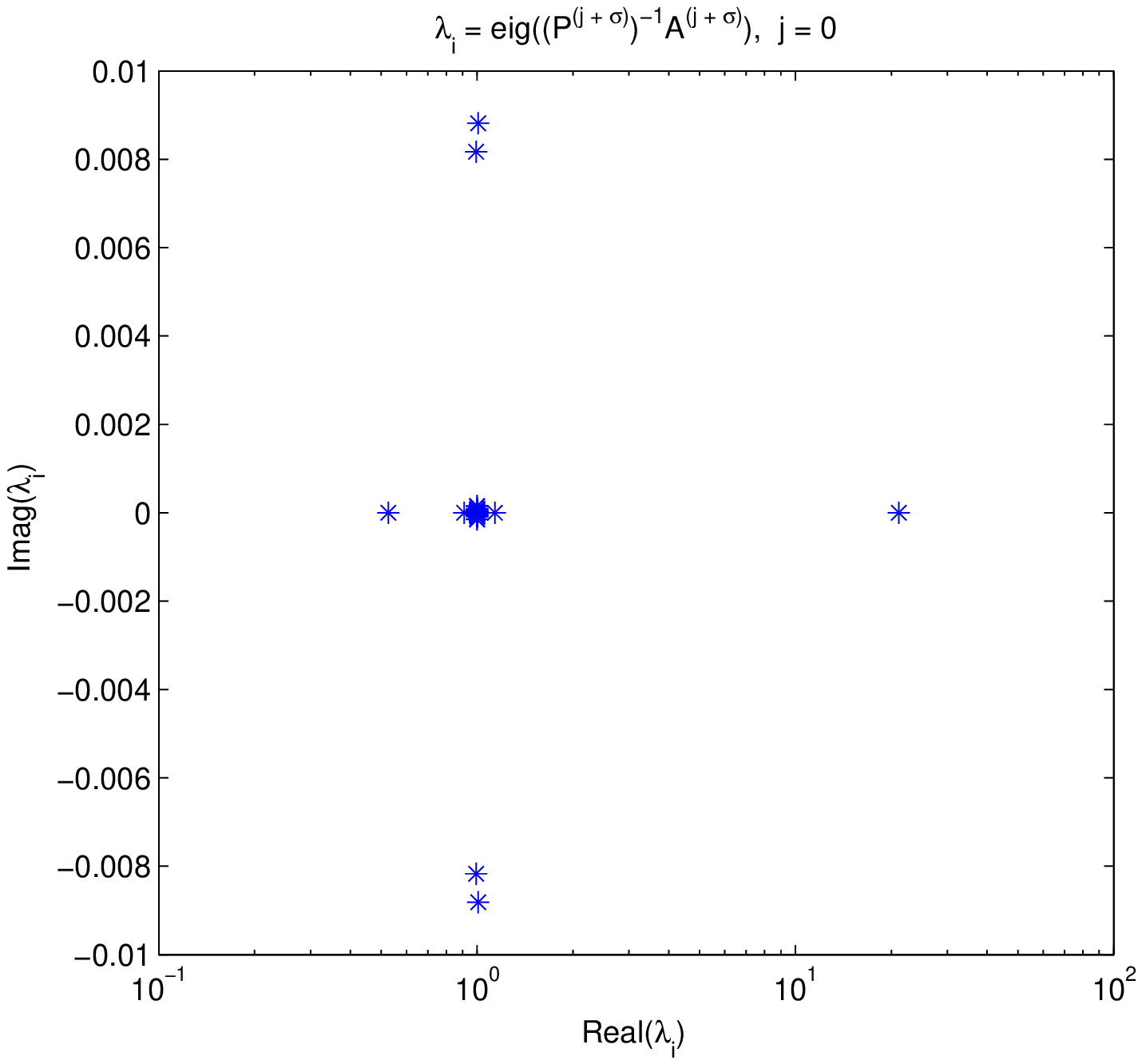}
\caption{{\small Spectrum of both original and preconditioned
matrices at the time level $j = 0$, respectively, when $N =
M = 128,~\alpha = 0.9$ and $\beta = 1.5$. Left: Original matrix;
Right: circulant preconditioned matrix.}}
\label{fig3}
\end{figure}

\begin{figure}[!htbp]
\centering
\includegraphics[width=3.1in,height=2.95in]{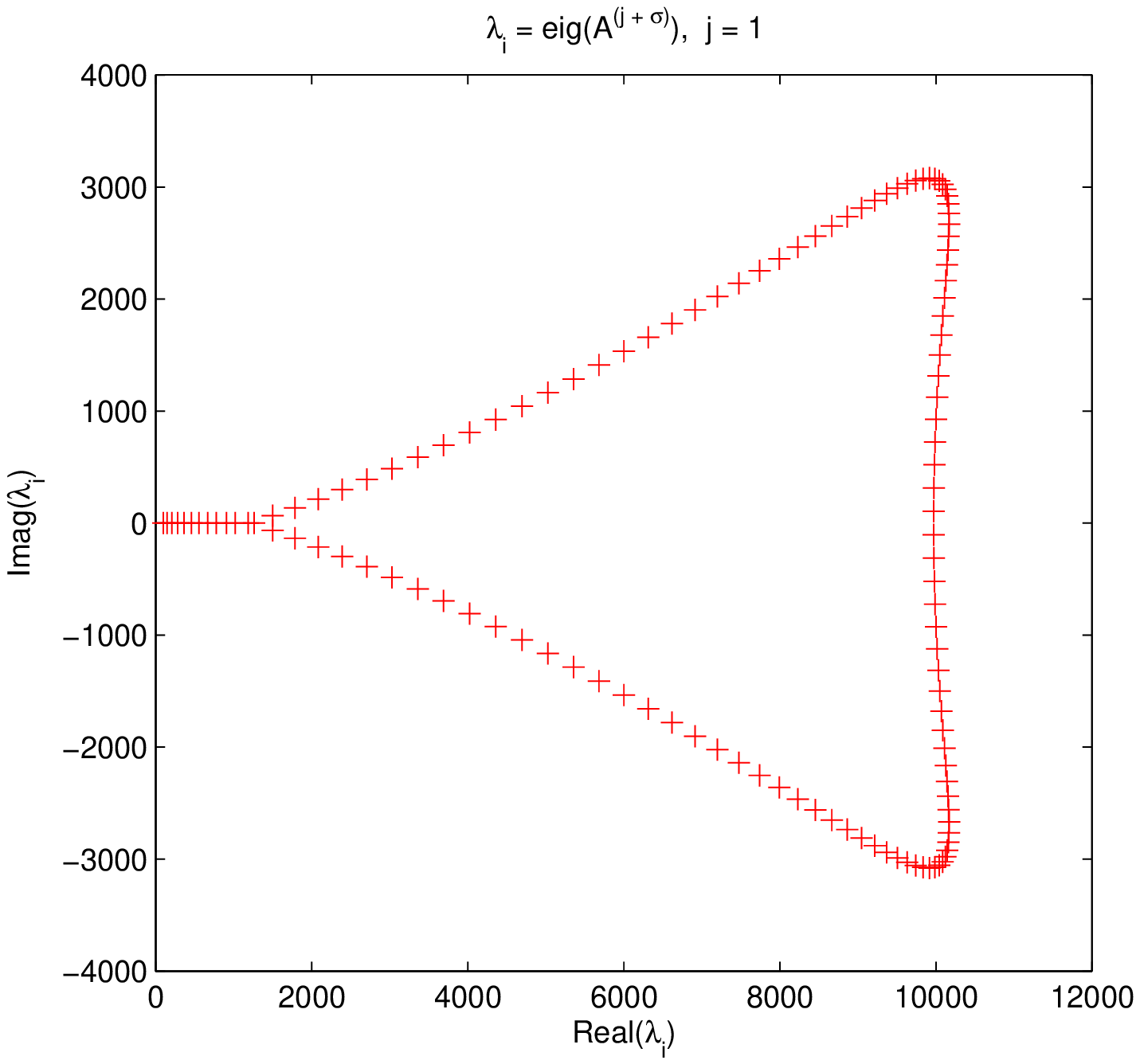}
\includegraphics[width=3.1in,height=2.95in]{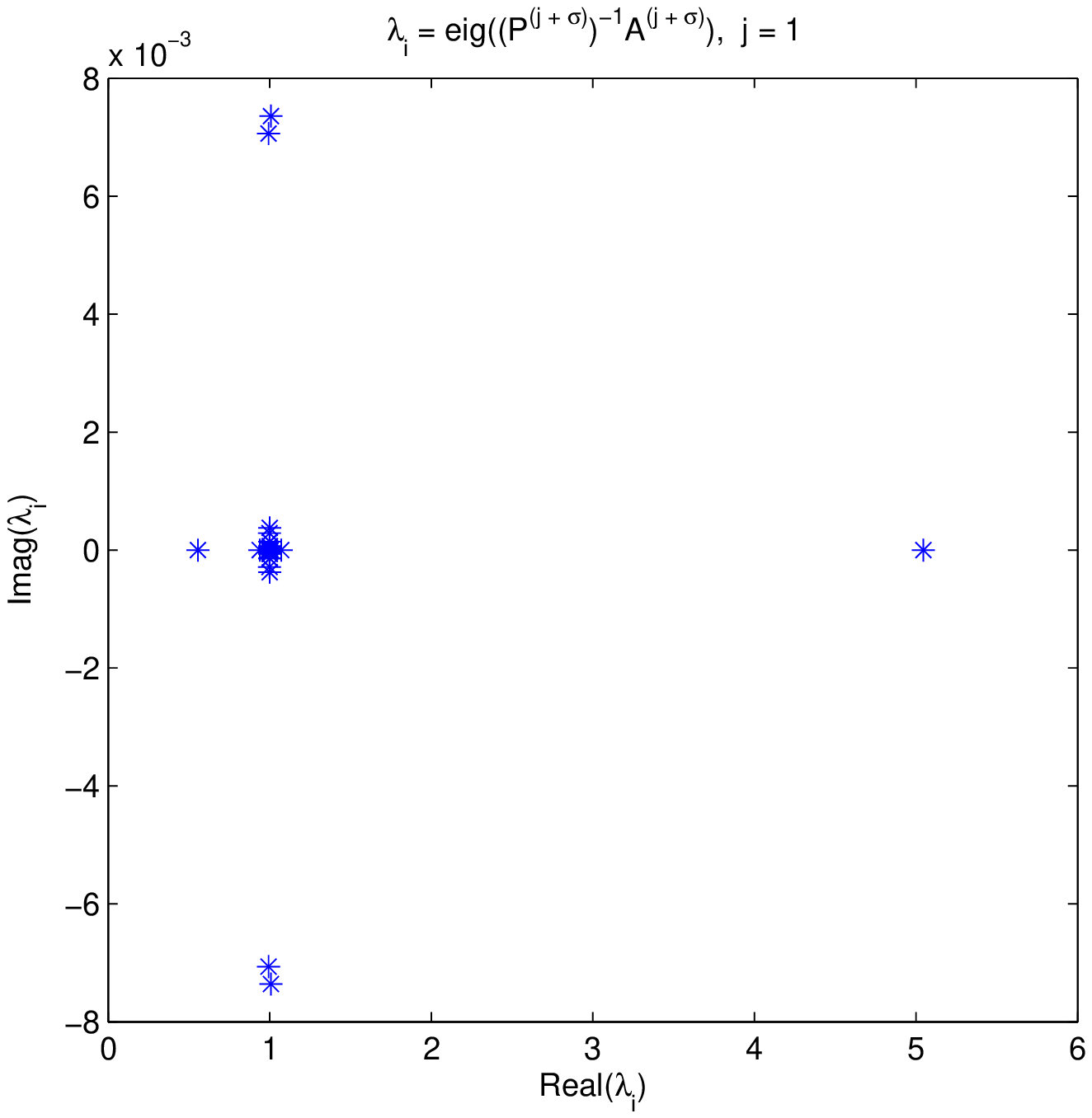}
\caption{{\small Spectrum of both original and preconditioned
matrices at the time level $j = 1$, respectively, when $N =
M = 128,~\alpha = 0.9$ and $\beta = 1.5$. Left: Original matrix;
Right: circulant preconditioned matrix.}}
\label{fig4}
\end{figure}

Again, for the case of variable time coefficients, several eigenvalue
plots about both original and preconditioned matrices are similarly
displayed in Figs. \ref{fig3}-\ref{fig4}. These two figures confirm that for
circulant preconditioning, the eigenvalues of preconditioned matrices are
clustered at 1, expect for few (about $6\sim 10$) outliers. The vast majority
of the eigenvalues are well separated away from 0. It may be mainly interpreted
as that in our implementation the number of iterations needed by PCGS with
circulant preconditioners almost ranges from 6 to 10. We validate the effectiveness
and robustness of the proposed circulant preconditioner from the perspective of
clustering spectrum.

\begin{table}[!htbp]\small\tabcolsep=5.3pt
\begin{center}
\caption{{\small CPU time in seconds for solving Example 2 with $\alpha = 0.9$ that
\texttt{Time1} is done by for Algorithm \ref{alg1} (MATLAB's backslash) and \texttt{Time2}
is done by Algorithm \ref{alg1} with PCGS solver.}}
\begin{tabular}{lrrcrrcrrc}
\hline &\multicolumn{3}{c}{$\beta = 1.3$} &\multicolumn{3}{c}{$\beta = 1.5$} &\multicolumn{3}
{c}{$\beta = 1.9$}\\
[-2pt]\cmidrule(l{0.7em}r{0.7em}){2-4} \cmidrule(l{0.7em}r{0.6em}){5-7}
\cmidrule(l{0.7em}r{0.7em}){8-10}\\[-11pt]
$h = \tau$ &\texttt{Time1} &\texttt{Time2} & $\texttt{Speed-up}$ &\texttt{Time1} &\texttt{Time2} &
$\texttt{Speed-up}$&\texttt{Time1} &\texttt{Time2} &$\texttt{Speed-up}$      \\
\hline
$2^{-5}$  &0.03   &0.06   &0.50  &0.03   &0.06  &0.50  &0.03   &0.06  &0.50  \\
$2^{-6}$  &0.09   &0.13   &0.69  &0.09   &0.14  &0.64  &0.09   &0.13  &0.69  \\
$2^{-7}$  &0.46   &0.32   &1.44  &0.47   &0.34  &1.38  &0.45   &0.34  &1.32  \\
$2^{-8}$  &2.65   &0.86   &3.08  &2.61   &0.91  &2.87  &2.59   &0.88  &2.94  \\
$2^{-9}$  &31.39  &3.70   &8.48  &31.63  &3.93  &8.05  &31.33  &3.74  &8.38  \\
$2^{-10}$ &316.34 &22.39  &14.13 &321.05 &22.96 &13.98 &329.89 &22.74 &14.51 \\
\hline
\end{tabular}
\label{tab8y}
\end{center}
\end{table}

\begin{table}[!htbp]\small\tabcolsep=5.3pt
\begin{center}
\caption{{\small CPU time in seconds for solving Example 2 with $\beta = 1.8$ that
\texttt{Time1} is done by for Algorithm \ref{alg1} (MATLAB's backslash) and \texttt{Time2}
is done by Algorithm \ref{alg1} with PCGS solver.}}
\begin{tabular}{rrrcrrcrrc}
\hline &\multicolumn{3}{c}{$\alpha = 0.5$} &\multicolumn{3}{c}{$\alpha = 0.9$} &\multicolumn{3}
{c}{$\alpha = 0.99$}\\
[-2pt]\cmidrule(l{0.7em}r{0.7em}){2-4} \cmidrule(l{0.7em}r{0.6em}){5-7}\cmidrule(l{0.7em}r{0.7em}){8-10}\\[-11pt]
$h = \tau$ &\texttt{Time1} &\texttt{Time2} & $\texttt{Speed-up}$ &\texttt{Time1} &\texttt{Time2} &
$\texttt{Speed-up}$&\texttt{Time1} &\texttt{Time2} &$\texttt{Speed-up}$      \\
\hline
$2^{-5}$  &0.03   &0.06   &0.50  &0.03   &0.06  &0.50  &0.03   &0.06  &0.50  \\
$2^{-6}$  &0.09   &0.13   &0.69  &0.09   &0.13  &0.69  &0.09   &0.13  &0.69  \\
$2^{-7}$  &0.45   &0.35   &1.28  &0.46   &0.34  &1.35  &0.45   &0.34  &1.32  \\
$2^{-8}$  &2.54   &0.98   &2.59  &2.55   &0.89  &2.87  &2.57   &0.90  &2.86  \\
$2^{-9}$  &31.56  &4.14   &7.62  &31.53  &3.92  &8.04  &31.87  &3.88  &8.21  \\
$2^{-10}$ &328.59 &23.46  &14.01 &328.73 &23.15 &14.20 &329.37 &22.83 &14.41 \\
\hline
\end{tabular}
\label{tab8z}
\end{center}
\end{table}

In Tables \ref{tab8y}-\ref{tab8z}, it verifies that the proposed fast direct
solver for different discretized problems takes much less CPU time elapsed as
$M$ and $N$ become large. When $M = N = 2^{10}$ and different discretized
parameters, the CPU time of Algorithm \ref{alg2} by PCGS with circulant
preconditioners is about 23 seconds, the speedup is almost 14 times. Meanwhile,
although \texttt{Time1} required by Algorithm \ref{alg2} with MATLAB's
backslash for small test problems ($M = N = 32,64$) than \texttt{Time2}
needed by Algorithm \ref{alg2} with using the PCGS method, our proposed method
is still more attractive in aspects of lower memory requirement. Compared to
Algorithm \ref{alg2} with MATLAB's backslash, it highlighted that in the whole
procedure the proposed solution technique does not require to store a series of
full matrices (e.g. coefficient matrices $A^{(j + \sigma)},~j = 0,1,\ldots,M - 1$) at all.
All in all, we can conclude that our proposed IDS with fast solution techniques
is still more promising than the IDS with common implementation.

\section{Conclusions}
\quad\
\label{sec5}
In this paper, the stability and convergence of an implicit difference schemes approximating the
time-space fractional convection-diffusion equation of a general form is studied. Sufficient
conditions for the unconditional stability of such difference schemes are obtained. For proving
the stability of a wide class of difference schemes approximating the time fractional diffusion
equation, it is simple enough to check the stability conditions obtained in this paper. Meanwhile,
the new difference schemes of the second approximation order in space and the second approximation
order in time for the TSFCDE with variable coefficients (in terms of $t$) are constructed as well.
The stability and convergence of these implicit schemes in the mesh $L_2$-norm with the rate equal to the
order of the approximation error are proved. The method can be easily adopted to other TSFCDEs
with other boundary conditions. Numerical tests completely confirming the obtained theoretical
results are carried out.

More significantly, with the aid of (\ref{eq3.1}), we can ameliorate  the calculation
skill by the implementation of reliable preconditioning iterative techniques, with only
computational cost and memory of $\mathcal{O}(N\log N)$ and $\mathcal{O}(N)$, respectively.
Extensive numerical results are reported to illustrate that the efficiency of the proposed
preconditioning methods. In future work, we will focus on the extension of the proposed
IDS for handling two/three-dimensional TSFCDEs with fast solution techniques subject to
various boundary value conditions. Meanwhile, we will also focus on the development of
other efficient preconditioners for accelerating the convergence of Krylov subspace solver
for the discretized Toeplitz systems; refer, e.g., to our recent work \cite{XMGTZH} for
this topic.

\section*{Acknowledgement}
\quad\
{\it We are grateful to Prof. Zhi-Zhong Sun and Dr. Zhao-Peng Hao for their
insightful discussions about the convergence analysis of the proposed implicit
difference scheme. The work of the first two authors is supported by 973 Program
(2013CB329404), NSFC (61370147, 11301057, and 61402082), the Fundamental Research
Funds for the Central Universities (ZYGX2013J106, ZYGX2013Z005, and ZYGX2014J084).
The work of the last author has been partially implemented with the financial
support of the Russian Presidential grant for young scientists MK-3360.2015.1.}


\begin{thebibliography}{99}

\bibitem{IPFDE}
I. Podlubny, Fractional Differential Equations, vol. 198 of Mathematics in Science
and Engineering, Academic Press, Inc., San Diego, CA, 1999.

\bibitem{SSAKO}
S.G. Samko, A.A. Kilbas, O.I. Marichev, Fractional Integrals and Derivatives: Theory
and Applications, Gordon and Breach Science Publishers, Yverdonn, 1993.

\bibitem{AAKHMS}
A.A. Kilbas, H.M. Srivastava, J.J. Trujillo, Theory and Applications of Fractional
Differential Equations, Elsevier, Amsterdam, 2006.

\bibitem{RMJKT}
R. Metzler, J. Klafter, The random walk's guide to anomalous diffusion: a fractional
dynamics approach, {\it Phys. Rep.}, 339 (2000), pp. 1-77.

\bibitem{AISGMZ}
A.I. Saichev, G.M. Zaslavsky, Fractional kinetic equations: solutions and applications,
{\it Chaos}, 7 (1997), pp. 753-764. Available online at \url{http://dx.doi.org/10.1063/1.166272}.

\bibitem{YPST}
Y. Povstenko, Space-time-fractional advection diffusion equation in a plane, in
{\em Advances in Modelling and Control of Non-integer-Order Systems}, K.J. Latawiec,
M. {\L}ukaniszyn, R. Stanis{\l}awski eds., Volume 320 of the series Lecture Notes in
Electrical Engineering, Springer, 2015, pp. 275-284.

\bibitem{DABSW}
D.A. Benson, S.W. Wheatcraft, M.M. Meerschaert, Application of a fractional advection-dispersion
equation, {\it Water Resour. Res.}, 36 (2000), pp. 1403-1412.

\bibitem{YZPF}
Y.Z. Povstenko, Fundamental solutions to time-fractional advection diffusion equation
in a case of two space variables, {\it Math. Probl. Eng.}, 2014 (2014), Article ID 705364,
7 pages. Available online at \url{http://dx.doi.org/10.1155/2014/705364}.

\bibitem{MMMCT}
M.M. Meerschaert, C. Tadjeran, Finite difference approximations for fractional
advection-dispersion flow equations, {\it J. Comput. Appl. Math.}, 172 (2004),
pp. 65-77.

\bibitem{LSPC}
L. Su, W. Wang, Q. Xu, Finite difference methods for fractional dispersion equations,
{\it Appl. Math. Comput.}, 216 (2010), pp. 3329-3334.

\bibitem{LSWWZ}
L. Su, W. Wang, Z. Yang, Finite difference approximations for the fractional
advection-diffusion equation, {\it Phys. Lett. A}, 373 (2009), pp. 4405-4408.

\bibitem{ESAS}
E. Sousa, Finite difference approximations for a fractional advection diffusion
problem, {\it J. Comput. Phys.}, 228 (2009), pp. 4038-4054.

\bibitem{LSWWHW}
L. Su, W. Wang, H. Wang, A characteristic difference method for the transient
fractional convection-diffusion equations, {\it Appl. Numer. Math.}, 61 (2011),
pp. 946-960.

\bibitem{KWHWA}
K. Wang, H. Wang, A fast characteristic finite difference method for fractional
advection-diffusion equations, {\it Adv. Water Resour.}, 34 (2011), pp. 810-816.

\bibitem{ZDVSLB}
Z. Deng, V. Singh, L. Bengtsson, Numerical solution of fractional advection-dispersion
equation, {\it J. Hydraul. Eng.}, 130 (2004), pp. 422-431.

\bibitem{ZDAXM}
Z. Ding, A. Xiao, M. Li, Weighted finite difference methods for a class of space
fractional partial differential equations with variable coefficients, {\it J. Comput.
Appl. Math.}, 233 (2010), pp. 1905-1914.

\bibitem{FLPZKB}
F. Liu, P. Zhuang, K. Burrage, Numerical methods and analysis for a class of
fractional advection-dispersion models, {\it Comput. Math. Appl.}, 64 (2012),
pp. 2990-3007.

\bibitem{SMAAR}
S. Momani, A.A. Rqayiq, D. Baleanu, A nonstandard finite difference scheme for
two-sided space-fractional partial differential equations, {\it Int. J. Bifurcat.
Chaos}, 22 (2012), 1250079, 5 pages. Available online at
\url{http://dx.doi.org/10.1142/S0218127412500794}.


\bibitem{MCWD}
M. Chen, W. Deng, A second-order numerical method for two-dimensional two-sided space
fractional convection diffusion equation, {\it Appl. Math. Model.}, 38 (2014), pp. 3244-3259.

\bibitem{WDMCE}
W. Deng, M. Chen, Efficient numerical algorithms for three-dimensional fractional partial
differential equations, {\it J. Comp. Math.}, 32 (2014), pp. 371-391.

\bibitem{WQSLL}
W. Qu, S.-L. Lei, S.-W. Vong, Circulant and skew-circulant splitting iteration for
fractional advection-diffusion equations, {\it Int. J. Comput. Math.}, 91 (2014),
pp. 2232-2242.

\bibitem{NJFLP}
N.J. Ford, K. Pal, Y. Yan, An algorithm for the numerical solution of two-sided
space-fractional partial differential equations, {\it Comput. Methods Appl.
Math.}, 15 (2015), pp. 497-514.

\bibitem{AHBDBA}
A.H. Bhrawy, D. Baleanu, A spectral Legendre-Gauss-Lobatto collocation method
for a space-fractional advection diffusion equations with variable coefficients,
{\it Rep. Math. Phys.}, 72 (2013), pp. 219-233.

\bibitem{AHBMAZ1}
A.H. Bhrawy, M.A. Zaky, A method based on the Jacobi tau approximation for solving
multi-term time-space fractional partial differential equations, {\it J. Comput.
Phys.}, 281 (2015), pp. 876-895.

\bibitem{HHTMFLS}
H. Hejazi, T. Moroney, F. Liu, Stability and convergence of a finite volume method
for the space fractional advection-dispersion equation, {\it J. Comput. Appl. Math.},
255 (2014), pp. 684-697.

\bibitem{WYTWD2}
W.Y. Tian, W. Deng, Y. Wu, Polynomial spectral collocation method for space fractional
advection-diffusion equation, {\it Numerical Methods for Partial Differential Equations},
30 (2014), pp. 514-535.

\bibitem{YLCXF}
Y. Lin, C. Xu, Finite difference/spectral approximations for the time-fractional diffusion
equation, J. Comput. Phys., 225 (2007), pp. 1533-1552.

\bibitem{MCAH}
M. Cui, A high-order compact exponential scheme for the fractional convection-diffusion
equation, {\it J. Comput. Appl. Math.}, 255 (2014), pp. 404-416.

\bibitem{MCCES}
M. Cui, Compact exponential scheme for the time fractional convection-diffusion
reaction equation with variable coefficients, {\it J. Comput. Phys.}, 280 (2015),
pp. 143-163.

\bibitem{AMMAC}
A. Mohebbi, M. Abbaszadeh, Compact finite difference scheme for the solution of
time fractional advection-dispersion equation, {\it Numer. Algorithms}, 63 (2013),
pp. 431-452.

\bibitem{SMAA}
S. Momani, An algorithm for solving the fractional convection-diffusion equation
with nonlinear source term, {\it Commun. Nonlinear Sci. Numer. Simul.}, 12 (2007),
pp. 1283-1290.

\bibitem{ZWSV}
Z. Wang, S. Vong, A high-order exponential ADI scheme for two dimensional time
fractional convection-diffusion equations, {\it Comput. Math. Appl.}, 68 (2014),
pp. 185-196.

\bibitem{ZJFWC}
Z.-J. Fu, W. Chen, H.-T. Yang, Boundary particle method for Laplace transformed time
fractional diffusion equations, {\it J. Comput. Phys.}, 235 (2013), pp. 52-66.

\bibitem{SZXFYH}
S. Zhai, X. Feng, Y. He, An unconditionally stable compact ADI method for
three-dimensional time-fractional convection-diffusion equation, {\it J. Comput.
Phys.}, 269 (2014), pp. 138-155.

\bibitem{PZYTG}
P. Zhuang, Y.T. Gu, F. Liu, I. Turner, P.K.D.V. Yarlagadda, Time-dependent fractional
advection-diffusion equations by an implicit MLS meshless method, {\it Int. J. Numer.
Meth. Eng.}, 88 (2011), pp. 1346-1362.

\bibitem{YMWAC}
Y.-M. Wang, A compact finite difference method for solving a class of time fractional
convection-subdiffusion equations, {\it BIT}, 55 (2015), pp. 1187-1217.

\bibitem{YZAFD}
Y. Zhang, A finite difference method for fractional partial differential equation,
{\it Appl. Math. Comput.}, 215 (2009), pp. 524-529.

\bibitem{YZFDA}
Y. Zhang, Finite difference approximations for space-time fractional partial
differential equation, {\it J. Numer. Math.}, 17 (2009), pp. 319-326.

\bibitem{YSWM}
Y. Shao, W. Ma, Finite difference approximations for the two-side space-time fractional
advection-diffusion equations, {\it J. Comput. Anal. Appl.}, 21 (2016), pp. 369-379.

\bibitem{PQXZA}
P. Qin, X. Zhang, A numerical method for the space-time fractional convection-diffusion
equation, {\it Math. Numer. Sin.}, 30 (2008), pp. 305-310. (in Chinese)

\bibitem{FLPZV}
F. Liu, P. Zhuang, V. Anh, I. Turner, K. Burrage, Stability and convergence of the difference
methods for the space-time fractional advection-diffusion equation, {\it Appl. Math. Comput.}, 191
(2007), pp. 12-20.

\bibitem{ZZXQJ}
Z. Zhao, X.-Q. Jin, M.M. Lin, Preconditioned iterative methods for space-time fractional
advection-diffusion equations, {\it J. Comput. Phys.}, 319 (2016), pp. 266-279.

\bibitem{SSFLVA}
S. Shen, F. Liu, V. Anh, Numerical approximations and solution techniques for
the space-time Riesz-Caputo fractional advection-diffusion equation, {\it Numer.
Algorithms}, 56 (2011), pp. 383-403.

\bibitem{MPMRE}
M. Parvizi, M.R. Eslahchi, M. Dehghan, Numerical solution of fractional
advection-diffusion equation with a nonlinear source term, {\it Numer.
Algorithms}, 68 (2015), pp. 601-629.

\bibitem{YCYWY}
Y. Chen, Y. Wu, Y. Cui, Z. Wang, D. Jin, Wavelet method for a class of fractional
convection-diffusion equation with variable coefficients, {\it J. Comput. Sci.}, 1
(2010), pp. 146-149.

\bibitem{SIMDS}
S. Irandoust-pakchin, M. Dehghan, S. Abdi-mazraeh, M. Lakestani, Numerical solution
for a class of fractional convection-diffusion equations using the flatlet oblique
multiwavelets, {\it J. Vib. Control}, 20 (2014), pp. 913-924.


\bibitem{AHBMAZ}
A.H. Bhrawy, M.A. Zaky, J.A. Tenreiro Machado, Efficient Legendre spectral tau
algorithm for solving the two-sided space-time Caputo fractional advection-dispersion
equation, {\it J. Vib. Control}, 22 (2016), pp. 2053-2068.

\bibitem{HHTMFL}
H. Hejazi, T. Moroney, F. Liu, A finite volume method for solving the two-sided
time-space fractional advection-dispersion equation, {\it Cent. Eur. J. Phys.},
11 (2013), pp. 1275-1283.

\bibitem{WJYLA}
W. Jiang, Y. Lin, Approximate solution of the fractional advection-dispersion
equation, {\it Comput. Phys. Commun.}, 181 (2010), pp. 557-561.

\bibitem{JWYCBL}
J. Wei, Y. Chen, B. Li, M. Yi, Numerical solution of space-time fractional
convection-diffusion equations with variable coefficients using Haar wavelets,
{\it Comput. Model. Eng. Sci. (CMES)}, 89 (2012), pp. 481-495.

\bibitem{MNIM}
M. Ng, Iterative Methods for Toeplitz Systems, Oxford University Press, New York,
2004.

\bibitem{FRLSWY}
F.-R. Lin, S.-W. Yang, X.-Q. Jin, Preconditioned iterative methods for fractional
diffusion equation, {\it J. Comput. Phys.}, 256 (2014), pp. 109-117.

\bibitem{SLLHWS}
S.-L. Lei, H.-W. Sun, A circulant preconditioner for fractional diffusion equations,
{\it J. Comput. Phys.}, 242 (2013), pp. 715-725.

\bibitem{IGAS}
I. Gohberg, A. Semencul, On the inversion of finite Toeplitz matrices and their continuous
analogues, {\it Matem. Issled.}, 7 (1972), pp. 201-223. (in Russian)

\bibitem{YSMHS}
Y. Saad, Iterative Methods for Sparse Linear Systems, 2nd edition, SIAM, Philadelphia,
USA, 2003.

\bibitem{4A2015}
A.A. Alikhanov, A new difference scheme for the time fractional diffusion equation,
{\it J. Comput. Phys.}, 280 (2015), pp. 424-438.

\bibitem{ZPHZZS}
Z.-P. Hao, Z.-Z. Sun, W.-R. Cao, A fourth-order approximation of fractional derivatives
with its applications, {\it J. Comput. Phys.}, 281 (2015), pp. 787-805.

\bibitem{SVPLX}
S. Vong, P. Lyu, X. Chen, S.-L. Lei, High order finite difference method for time-space
fractional differential equations with Caputo and Riemann-Liouville derivatives, {\it Numer.
Algorithms}, 72 (2016), pp. 195-210.

\bibitem{XMGTZH}
X.-M. Gu, T.-Z. Huang, H.-B. Li, L. Li, W.-H. Luo, On $k$-step CSCS-based polynomial
preconditioners for Toeplitz linear systems with application to fractional diffusion
equations, {\it Appl. Math. Lett.}, 42 (2015), pp. 53-58.

\bibitem{GHGCFV}
{\AA}. Bj\"{o}rck, Direct Methods for Linear Systems, in \textit{Numerical Methods in
Matrix Computations}, Volume 59 of the series Texts in Applied Mathematics, Springer,
Switzerland, 2014, pp. 1-209.
\end{thebibliography}
\end{document}